\documentclass[11pt]{amsart}
\usepackage[all]{xy}

\usepackage[doi=false,isbn=false,url=false
]{biblatex}
\theoremstyle{plain}
\newtheorem*{Thm3.1}{Theorem~3.1}
\newtheorem*{Thm4.1}{Theorem~4.1}
\newtheorem*{Thm4.2}{Theorem~4.2}

\usepackage[a4paper, total=
{6.23in,9in}]{geometry}

\title{Contact homology of contact manifolds and its applications}

\author{Fr\'ed\'eric Bourgeois} 
 \address{Universit\'e Paris-Saclay, CNRS, Laboratoire de Math\'ematiques d'Orsay, 
91405 Orsay, France} 
\urladdr{https://www.imo.universite-paris-saclay.fr/{\raisebox{0.5ex}{\texttildelow}}frederic.bourgeois/}
\email{frederic.bourgeois@universite-paris-saclay.fr}

\thanks{The author's research activity is partially supported by the ANR project COSY (21-CE40-0002).}

\numberwithin{equation}{section}

\begin{document}

\begin{abstract}
This is a survey of contact homology and its applications to the study of contact manifolds. It is a small tribute to Yasha Eliashberg's huge 
generosity with his countless explanations of his deep mathematical insights all along his career. It is also the author's wishful thinking that this
text could be useful to students and young mathematicians for learning about some of the holomorphic curves based invariants in 
contact geometry.
\end{abstract} 

\maketitle

\section*{1. Introduction}

In order to describe the context in which contact homology was born, it is useful to go back to the origins of a whole collection of symplectic and contact 
invariants. Holomorphic curves were introduced in symplectic geometry by 
Gromov
\cite{mr:809718},
in the case of closed symplectic manifolds $(W, \omega)$. 
These are maps $u : \Sigma \to W$ defined on a Riemann surface $\Sigma$ with
complex structure $j$ satisfying the Cauchy--Riemann equation
\[
du \circ j = J \circ du
\]
for a compatible almost complex structure $J$ on $W$, i.e., an endomorphism $J$ of $TW$ satisfying $J^2 = -I$ as well as the 
compatibility conditions 
\begin{itemize}
\item $\omega(J \cdot\,, J \cdot\,) = \omega(\,\cdot\,, \cdot\,)$, 
\item $\omega(v, Jv) > 0$ for all nonzero tangent vectors $v$ to $W$.
\end{itemize}
Sometimes the first condition is omitted, and we then say that our almost complex structure $J$ is tame.

Originally, such maps $u$ were referred to as pseudoholomorphic, because the compatible almost complex structure $J$ is generally not integrable,
i.e., one cannot find complex local coordinates on $W$ so that $J$ corresponds to the multiplication by $i$ in these coordinates. As mathematicians became 
used to these objects with time,  the name $J$-holomorphic became fashionable. Nowadays, one can simply refer to these maps as ``holomorphic'' if the 
context is clear. This last simplification may even be 
due to Yasha Eliashberg's pleas during 
various conferences and lectures, where he has insisted 
with an evident sense of humor that the prefixes ``pseudo-'' or ``$J$-'' suggest that these maps are not as good as actual holomorphic 
ones, while in fact they are perfectly fine, and that this terminology is ``degrading our field''. He even suggested 
calling these maps ``Gromomorphic'' 
but for some reason it seems that this name was never used in the literature. 

In a nutshell, the fact that moduli spaces of $J$-holomorphic curves can lead to symplectic invariants mainly rests on two key facts. 
First, the set of compatible almost complex structures on a symplectic manifold
$(W, \omega)$ is contractible. 
In particular, any symplectic manifold can be equipped with such a compatible
almost complex structure.  Moreover, different choices of compatible almost
complex structures are always isotopic, and different isotopies between fixed
compatible almost complex structures are themselves isotopic, and so on.

Second, the moduli spaces of $J$-holomorphic curves can be compactified by
adding nodal curves, i.e., Riemann surfaces with pairs of points 
identified to each other and called \textit{nodes}, as well as $J$-holomorphic maps from these Riemann surfaces to $W$ that are compatible with these 
identifications. This is the content of the celebrated Gromov compactness theorem in symplectic geometry. In particular, intersection numbers in 
these compactified moduli spaces are well defined and can be used in the construction of symplectic invariants. In some cases, these intersections 
numbers are the symplectic invariants themselves; this is the case for Gromov--Witten invariants. In other situations, these 
intersection numbers 
are used as coefficients in a differential of a chain complex, which is defined in such a way that its homology is a symplectic invariant; this is the 
case of Floer homology. For the latter type of theory, nodal curves are replaced with broken trajectories in the compactification of moduli spaces. 
The simplest type of broken trajectory consists of two cylinders, such that the end of the first one coincides with the beginning of the second one.

Some might say that Floer homology is not strictly speaking an invariant based on holomorphic curves, since the corresponding differential counts 
solutions of the so-called Floer equation, which is an inhomogeneous Cauchy--Riemann equation, involving a Hamiltonian function in its right-hand 
side. But as explained in 
\cite[paragraph~1.4.C$'$]{mr:809718},
the graph of a map satisfying such an inhomogeneous Cauchy--Riemann equation is 
$\widetilde{J}$-holomorphic for a suitable choice of almost complex structure $\widetilde{J}$ on the target space of this graph. Yasha indeed 
insists that it is very useful to read the whole of Gromov's original article 
\cite{mr:809718}
as it contains many fruitful ideas. More specifically, this 
interpretation of Floer trajectories as holomorphic maps is explained in 
\cite[Section~2.2]{mr:2026549}
using the notion of a stable Hamiltonian structure. 
All this is a beautiful illustration of one of Yasha's mottos: ``Look at
graph'' which is 
fundamental in an impressive number of situations within contact 
and symplectic topology.

When $(W, \omega)$ is a symplectic manifold with boundary, satisfying some type of convexity property near its boundary, holomorphic curves 
obey a maximum principle, 
hence they are confined in a compact region of the target manifold $W$, and the corresponding moduli spaces are 
still compact. This is, for example, the case of symplectic manifolds $(W, \omega)$ with contact-type boundary $(\partial W = M, \xi)$, for which 
there exists near $\partial W$ a Liouville vector field $v$, i.e., $\mathcal{L}_v \omega = \omega$, such that $v$ is transverse to $\partial W$ and 
pointing outside $W$. In that case, the kernel of the $1$-form $\alpha = \imath(v) \omega$ is a contact structure $\xi$. In other words, the 
$1$-form $\alpha$ satisfies the property that $\alpha \wedge (d\alpha)^{n-1}$
does not vanish wherever it is defined, and is 
called a 
\textit{contact form}. 
Different choices $v_1$ and $v_2$ of a vector field as above lead to different contact structures, but the collection of all convex combinations of 
$v_1$ and $v_2$ leads to an isotopy between the corresponding contact
structures, and by 
Gray's stability theorem for closed contact manifolds, 
these contact structures are actually diffeomorphic.

Some beautiful instances of successful applications of holomorphic curves in this type of symplectic manifolds $(W, \omega)$ are described 
by Yasha in 
\cite{mr:1171908}.
There, holomorphic maps are defined over the disc and its boundary must be mapped to a given submanifold of 
$(W, \omega)$ satisfying appropriate conditions. The development of these techniques were followed by an impressive number of further developments. 
Those are described in another chapter of this volume \cite{cm:geiges.2024}. 
A more recent theory relying on holomorphic curves in such $(W, \omega)$ is symplectic homology, which is a version of Floer homology using a 
Hamiltonian function and a compatible almost complex structure having an appropriate asymptotic behavior near $\partial W$.

In a symplectic manifold $(W, \omega)$ with contact type boundary, any sufficiently small open neighborhood of $\partial W = M$ can be symplectically 
embedded in the symplectization of $(M, \xi)$. The latter is the manifold $\mathbb{R} \times M$ equipped with the symplectic form $d(e^t \alpha)$, where $t$ 
is a global coordinate on the $\mathbb{R}$ factor and $\alpha$ is a contact form for $\xi$. Note that this symplectic manifold does not depend on the choice of $\alpha$.
Yasha used to explain 
to his students that, for some time, symplectizations were considered by experts
to be bad manifolds to work with holomorphic curves. 
This is because the lower end of a symplectization (with $t$ very small) is concave, as the Liouville vector field $\frac{\partial}{\partial t}$ is pointing inside the manifold, instead 
of being convex. Therefore, holomorphic curves are allowed to sink deeply into this lower end and visit an unbounded region of the manifold, so that the 
corresponding moduli spaces will not be compact. This all changed when 
Hofer
\cite{mr:1244912}
understood the behavior of finite energy holomorphic curves 
at infinity of a symplectization, when equipped with a compatible almost complex structure with an appropriate behavior outside a compact set. This will be
detailed below, but for now let us just say that near each end of a symplectization, such holomorphic curves are asymptotic to cylinders over periodic orbits 
of a distinguished vector field in the contact manifold.

Symplectizations are often schematically represented as cones, since the symplectic volume tends to zero as $t \to -\infty$. But then, when drawing 
a holomorphic curve sinking into the negative end of a symplectization, it gets impossible to see what happens to it as the picture becomes 
arbitrarily small.
Maybe this partially explains why the behavior of holomorphic curves
symplectizations was not understood before Hofer's breakthrough work. 
In a public address during the first Yashafest conference at Stanford in 2007, 
Alexander Givental
astutely observed with a lot of humor that Yasha's 
vision of mathematics is so great because he is very small. What this really means, he explained, is that Yasha imagines mathematical objects 
as being much larger than himself, so that he can visualize all 
their details and develop a good intuition. This is obvious for anyone who 
has had the 
chance to discuss mathematics with Yasha on a sidewalk 
terrace, where it would be typical for him to describe a huge holomorphic
cylinder, say,
coming from the other side of the street and illustrate with it some new idea that had crossed his mind. Coming back to symplectizations, it is a better 
idea to draw them as cylinders, so that one can see more clearly what happens to holomorphic curves at their negative end.

Shortly after the work of Hofer, Yasha and 
he envisioned how to adapt ideas from Floer homology to symplectizations: this was the birth of contact 
homology. Yasha described contact homology in the proceedings of his invited lecture at the ICM 1998 in Berlin 
\cite{mr:1648083}.
As we shall see below, 
this theory was geometrically and algebraically more complex than Floer theory, but this did not stop them from developing an even more general theory. 
With 
Givental,
they introduced Symplectic Field Theory (SFT) 
\cite{mr:1826267},
which can be thought of as the general framework 
for using holomorphic curves 
of arbitrary topology in symplectic manifolds with convex and concave ends. 
SFT is also described in an article by Yasha in the proceedings of the ICM 2006 
in Madrid 
\cite{mr:2334192}.
This very general theory is the subject of  another chapter in this volume
\cite{cm:hind.seigel.2024}. 

One 
of the more general features of 
full SFT in comparison with contact homology is that a pair of holomorphic curves seen as a broken configuration 
may be glued along an arbitrary number of pairs of cylindrical ends, instead of
along a single pair of ends. This additional complexity plays an important 
role in the description of the relevant compactified moduli spaces, as well as in the algebraic formalism of the theory. During the second SFT Workshop in 
Leipzig in 2006, Yasha gave many lectures on SFT formalism. In order to illustrate this more sophisticated gluing for holomorphic curves in this 
context, he repeatedly made the following gesture. He placed one hand with its fingers pointing up below the other hand with its fingers pointing down, 
so that matching fingertips would figure pairs of cylindrical ends that could potentially be glued to each other, and in order to convey the fact that gluing was 
optional for each pair, he wiggled all his fingers with his hands in this position. In order to answer most questions from audience members puzzled by 
SFT, Yasha would repeat the same gesture. This eventually prompted a desperate cry by 
Viktor Ginzburg
from the back of the 
amphitheater: 
``Yasha, define \textit{something}!'' As a 
testimony to Yasha's generosity, a special session was added to the program 
at Yasha's request so that he could explain in more detail these 
puzzling aspects of SFT.  This session lasted until so late at night that the gates outside the building were locked, and all participants 
(including Yasha himself) had to climb over the wall in order to leave the campus.

Since the author cannot wiggle his fingers in front of the 
reader, a few things
will be defined in Sections~2 
and 3, 
concerning holomorphic curves in symplectizations and the contact homology
algebra, respectively, in the hope that readers will be able to comprehend
the applications of contact homology 
explained in Section~4. 

\section*{2. Holomorphic curves in symplectizations}

In this section, we describe geometric constructions in the symplectization of a contact manifold, leading to the definition of a suitable class of 
holomorphic curves for our purposes. We then outline the analysis required to construct the moduli spaces of such holomorphic curves with
the necessary geometric properties in view of the definition of homological invariants.

\subsection*{\textup{2.1.} Closed Reeb orbits}
\label{sec:orbits}

Let $(M, \xi)$ be a closed contact manifold of dimension $2n-1$. To any contact form $\alpha$ for $\xi$, one can associate its Reeb vector field
$R_\alpha$, characterized by the property that it spans the 
one-dimensional kernel of $d\alpha$: $\imath(R_\alpha) \,d\alpha = 0$, and by the 
normalization condition $\alpha(R_\alpha) = 1$. Although this vector field strongly depends on the choice of $\alpha$, the collection of Reeb fields for 
a given contact structure 
shares important dynamical properties. 

Of particular interest are the periodic orbits of such a vector field. Let 
$\gamma : \mathbb{R}/T\mathbb{Z} \to M$ be a closed Reeb orbit of period $T$: $\dot{\gamma}(t) = R_\alpha(\gamma(t))$. In other words, if we denote the flow 
of $R_\alpha$ by $(\varphi^\alpha_t)_{t \in \mathbb{R}}$, we have $\varphi^\alpha_T(\gamma(0)) = \gamma(0)$. Since $\mathcal{L}_{R_\alpha} \alpha = 0$, 
this flow preserves the contact structure, and we can restrict $d\varphi^\alpha_T$ to $\xi_{\gamma(0)}$. The closed Reeb orbit $\gamma$ is said to be 
nondegenerate if this restriction does not have $1$ as an eigenvalue. It can be shown that for a generic choice of $\alpha$, all closed Reeb orbits are 
nondegenerate; in this case, the contact form $\alpha$ is said to be nondegenerate. We denote by $\mathcal{P}_\alpha$ the collection of all 
periodic orbits of the Reeb field $R_\alpha$ modulo reparametrization shift, and for any $T > 0$ we denote by $\mathcal{P}^{\le T}_\alpha$ the 
subset of closed orbits with period less than or equal to $T$. For a nondegenerate contact form $\alpha$, the set $\mathcal{P}^{\le T}_\alpha$ 
is finite for all $T > 0$, and the set $\mathcal{P}_\alpha$ is at most countable. 

To any nondegenerate closed Reeb orbit $\gamma$, given a symplectic 
trivialization $\mathcal{F}$ of $\gamma^*\xi$, one can associate its Conley--Zehnder index $\mu_{CZ}(\gamma,
\mathcal{F}) \in \mathbb{Z}$. Intuitively, this 
index measures the twisting of the Reeb flow around the closed orbit $\gamma$, relative to the chosen trivialization $\mathcal{F}$. Any other symplectic 
trivialization $\mathcal{F}'$ can be obtained by acting on $\mathcal{F}$ with a loop of symplectic matrices $P :
\mathbb{R}/T\mathbb{Z} \to \mathrm{Sp}(2n-2)$. A natural
identification of $\pi_1(\mathrm{Sp}(2n-2))$ with $\mathbb{Z}$ is obtained by retracting a loop in $\mathrm{Sp}(2n-2)$ to $\mathrm{U}(n-1)$, then taking its 
complex determinant with values 
in $S^1 \subset \mathbb{C}$ and considering the degree of the resulting map from $S^1$ to itself. In particular, to the above loop $P$ one can associate in this 
way its Maslov index $\mu(P) \in \mathbb{Z}$. It then turns out that 
\[
\mu_{CZ}(\gamma, \mathcal{F}') = \mu_{CZ}(\gamma, \mathcal{F}) + 2 \mu(P).
\]
It follows that the Conley--Zehnder index depends only on the homotopy class of the chosen trivialization, and that its parity does not depend on it at all.

Note that there is a consistent way of obtaining symplectic trivializations for all orbits in $\mathcal{P}_\alpha$ from a minimal number of choices
(see for example 
\cite[Section~2.1]{mr:2555933}),
but we shall not use this here in order to keep the exposition more elementary. Instead, let us choose
once and for all a symplectic trivialization for each $\gamma \in
\mathcal{P}_\alpha$, and denote its Conley--Zehnder index with respect to this
chosen trivialization by $\mu_{CZ}(\gamma)$.

Let $\gamma$ be a closed Reeb orbit with minimal period $T$, i.e., the map $\gamma : [0, T) \to M$ is injective. We can consider its multiples 
\[
\gamma^{(m)} : \mathbb{R}/mT\mathbb{Z} \to M : t \mapsto \gamma(t \textrm{ mod }
T),
\]
for all integers $m \ge 1$. We say that a closed Reeb orbit is bad if it
is of the form $\gamma^{(m)}$ and if $\mu_{CZ}(\gamma^{(m)})$ and $\mu_{CZ}(\gamma^{(1)})$ have different parities; otherwise, we say that
it is good. Note that bad orbits are necessarily of even multiplicity $m$, and that if $\gamma^{(m)}$ is bad for some $m$, then all orbits of the form
$\gamma^{(2k)}$ are bad for all $k$. We define by $\mathcal{P}^{\rm g}_\alpha$ the subset of good closed Reeb orbits for the Reeb vector 
field corresponding to the contact form $\alpha$.

\subsection*{\textup{2.2.} Holomorphic curves}

Consider the symplectization $(\mathbb{R} \times M, d(e^t \alpha))$ of our contact manifold, where $t$ is a global
coordinate on $\mathbb{R}$. Note that the 
symplectomorphism class of this manifold does not depend on the choice of the contact form $\alpha$ for $\xi$. Choose an almost complex 
structure $J$ on this symplectization, i.e., an endomorphism of the tangent bundle of $\mathbb{R} \times M$ such that $J^2 = -I$, and satisfying the
following properties:
\begin{enumerate}
\item $J$ preserves $\xi$ and $J |_\xi$ is compatible with the symplectic form $d\alpha$, in the sense that 
	\[
		d\alpha(Jv,Jw) =d\alpha(v,w)
	\]
for all tangent vectors $v, w$ in $\xi$, and $d\alpha(v,Jv) > 0$ for all nonzero tangent vectors $v$ in $\xi$; 
\item $J \frac{\partial}{\partial t} = R_\alpha$;
\item $J$ is invariant by translation in the $\mathbb{R}$-direction.
\end{enumerate}
The first two conditions imply that $J$ is compatible with the symplectic structure $\omega = d(e^t \alpha)$.  
Like on closed symplectic manifolds, the space $\mathcal{J}(M,\alpha)$ of such almost complex 
structures is contractible.

Consider the Riemann sphere $\mathbb{C}\mathrm{P}^1$ equipped with its natural
complex structure $j$, as well as distinct points 
(called \textit{punctures}) 
$x_1, \dots, x_q \in \mathbb{C}\mathrm{P}^1$ with $q \ge 0$.
Given a map 
\[
u : \Sigma_u = \mathbb{C}\mathrm{P}^1 \setminus \{ x_1, \dots, x_q \} \to \mathbb{R} \times M
\]
defined on the punctured Riemann surface $\Sigma_u$, 
we denote by $u_\mathbb{R}$ and $u_M$ its components in
$\mathbb{R}$ and in $M$, respectively. Such a map is said to be $J$-holomorphic if $du \circ j = J \circ du$. The Hofer energy of $u$ is defined as
\[
E(u) = \sup_{\phi \in \mathcal{C}} \int_{\Sigma_u} u^*d(\phi \alpha),
\]
where $\mathcal{C} = \{ \phi \in
C^\infty(\mathbb{R}, [0,1] \ | \ \phi'(t) \ge 0 \}$.
This is a substitute for the symplectic area with respect to the symplectic form $d(e^t \alpha)$, which would be infinite if $u$ 
approached 
the positive end of the symplectization. Another interesting quantity is the so-called area of $u$, defined by
\[
A(u) = \int_{\Sigma_u} u^* \,d\alpha.
\]
Note that since $J$ is compatible with $d(e^t \alpha)$ and $J |_\xi$ is compatible with $d\alpha$,
the Hofer energy and the area are nonnegative: $E(u) \ge 0$ and $A(u) \ge 0$ for any holomorphic map $u$. Moreover, if $E(u)=0$, 
then $u$ is a constant map, and if $A(u)=0$, then the image of $u_M$ is contained in a Reeb trajectory.

It turns out that any proper holomorphic map $u$ with finite Hofer energy, i.e., $E(u) < \infty$, has a precise behavior near each of its punctures.
If $p$ denotes a puncture of $\Sigma_u$, the module and argument of a local complex coordinate for
$\mathbb{C}\mathrm{P}^1$ centered on $p$ determine
polar coordinates $(\rho, \theta)$ near $p$.  Then we have 
\[
\lim_{\rho \to 0} u_\mathbb{R}(\rho, \theta) = \pm \infty
\quad\text{ and }\quad
\lim_{\rho \to 0} u_M(\rho, \theta) = \gamma(\mp T\theta/2\pi),
\]
for some $\gamma \in \mathcal{P}_\alpha$ of period $T$. In this case, we say 
that the map $u$ is asymptotic to $\gamma$ at $\pm \infty$ and that the corresponding puncture $p$ is positive (resp. negative).

The prescribed asymptotic behavior of such maps near the $q$ punctures is specified by a list $\Gamma$ of $q$ periodic orbits in 
$\mathcal{P}_\alpha$. By convention, the asymptotes corresponding to positive punctures are listed first, and are separated by a 
semicolon from the asymptotes corresponding to negative punctures that are listed afterwards.
Two holomorphic maps $u$ and $u'$ with the same asymptotic behavior $\Gamma$ are said to be equivalent, 
which will be denoted by $u \sim u'$, if there exists a biholomorphism $h : \Sigma_u \to \Sigma_{u'}$ such that $u' \circ h = u$. The equivalence
class of a holomorphic map $u$ will be denoted by $[u]$.

The set of equivalence classes of holomorphic maps with a given asymptotic behavior $\Gamma$ will be denoted by 
$\widetilde{\mathcal{M}}(\Gamma)$. As this set actually parametrizes the solutions of partial differential equation modulo reparametrization, 
it deserves to 
be called a moduli space. Elements of this moduli space can be given a homology class in $H_2(M,
\mathbb{Z})$ after some appropriate 
choices have been fixed for all closed orbits in $\mathcal{P}_\alpha$. This moduli space can then be decomposed as a disjoint union of smaller 
moduli spaces according to these homology classes (see, for example,
\cite[Section~2.2]{mr:2555933}).
Again, we shall not make such choices here in order 
to keep the exposition more elementary.

Since the almost complex structure $J$ was chosen to be $\mathbb{R}$-invariant on the symplectization, the above moduli spaces are equipped with 
an $\mathbb{R}$-action: for any $\tau \in \mathbb{R}$, we set $\tau \cdot [u_\mathbb{R}, u_M] = [u_\mathbb{R} +
\tau, u_M]$ for each equivalence class $[u] = [u_\mathbb{R}, u_M]$ 
in some of these moduli spaces. The fixed points of this $\mathbb{R}$-action are easily identified as the vertical cylinders over periodic orbits in 
$\mathcal{P}_\alpha$. This corresponds to a very particular choice of asymptotic behavior $\Gamma = (\gamma ; \gamma)$, with one positive 
and one negative puncture each corresponding to the same closed orbit $\gamma$ of period $T$. In that case, it follows from Stokes theorem 
that for any $[u] \in \widetilde{\mathcal{M}}(\Gamma)$, we have $A(u) = T - T = 0$ so that the above $\mathbb{R}$-action is trivial. Except in this very 
special case of moduli spaces consisting of trivial cylinders, this $\mathbb{R}$-action is therefore free. We then consider the corresponding quotient 
$\widetilde{\mathcal{M}}(\Gamma)/\mathbb{R} = \mathcal{M}(\Gamma)$. In the case of trivial cylinders, we simply set 
$\mathcal{M}(\gamma;\gamma) = \widetilde{\mathcal{M}}(\gamma;\gamma)$.

Before turning to the properties of these moduli spaces $\mathcal{M}(\Gamma)$, let us specify the particular type of holomorphic curves which
will be of interest in order to define contact homology. As sketched in
Section~1, 
contact homology can be 
thought of as
the simplest theory for holomorphic curves adapting the ideas from Floer homology to symplectizations. Our aim is therefore to limit ourselves
to the simplest possible types of holomorphic curves. Note that in the above discussion, we already restricted ourselves to rational curves.
A first naive choice would be to further restrict ourselves to holomorphic cylinders, like in Floer homology where only Floer trajectories parametrized 
by $\mathbb{R} \times S^1$ are counted in order to define the differential. In the case of symplectizations, this would correspond to holomorphic curves with 
one positive and one negative 
puncture. However, this simple approach does not work in general, for the following reason. In order to construct
invariants out of moduli spaces, it is essential to consider and understand their compactifications. In our situation, we have to understand the 
compactification $\overline{\mathcal{M}}(\Gamma)$ of our moduli spaces $\mathcal{M}(\Gamma)$, using compactness results from 
\cite{mr:2026549}.
Roughly speaking, this compactification is obtained by adding to $\mathcal{M}(\Gamma)$ the so-called holomorphic buildings, consisting of 
several (but finitely many) levels. Each level contains one or more holomorphic curves in a symplectization, so that its positive asymptotes coincide 
with the negative asymptotes from the level directly above it, and its negative asymptotes coincide with the positive asymptotes from the level 
directly 
below it. The topmost level has positive asymptotes $\Gamma^+$ and the bottommost level has negative asymptotes $\Gamma^-$, such that
$\Gamma = (\Gamma^+ ; \Gamma^-)$. In particular, the compactification $\overline{\mathcal{M}}(\Gamma)$ of a moduli space $\mathcal{M}
(\Gamma)$ consisting of holomorphic cylinders could contain holomorphic
buildings with two levels, the top level containing a pair of pants, i.e.,
a holomorphic curve with one positive and two negative punctures, and the bottom level containing two holomorphic curves: a holomorphic 
cylinder and a holomorphic plane, i.e., a holomorphic curve with just one positive puncture. Note that this is just the simplest possible 
situation involving more general holomorphic curves than simply holomorphic cylinders, but we shall see later that it is of particular 
importance. For general contact manifolds, such holomorphic buildings can occur, so that it is not possible to work with compact moduli spaces
consisting of holomorphic buildings involving holomorphic cylinders only. Instead, one has to find a more general class of holomorphic curves
so that the corresponding compactified moduli spaces will solely consist of holomorphic buildings involving holomorphic curves in the same class.
If this class of holomorphic curves contains holomorphic cylinders, then it must
also contain tree-like holomorphic curves, i.e., holomorphic 
curves with exactly one positive puncture, but with an arbitrary number of negative punctures, as shown by the above example with more 
than one holomorphic plane in the bottom level. On the other hand, the compactification of moduli spaces of tree-like holomorphic curves cannot
involve holomorphic curves with more than one positive puncture. If this were the case, the corresponding holomorphic building would need to
contain at least one holomorphic curve $u$ without positive puncture but with at least one negative puncture. But such a holomorphic map $u$ 
cannot exist, because by Stokes theorem, the area $A(u)$ of this map would be negative, which is forbidden.

The above discussion explains why, for the purpose of defining contact homology, we shall require that our holomorphic maps $u$ have exactly 
one positive puncture and an arbitrary number $r \ge 0$ of negative punctures. Their asymptotic behavior will be of the form $\Gamma = 
(\gamma^+ ; \gamma^-_1, \dots, \gamma^-_r)$, where $\gamma^+ \in \mathcal{P}_\alpha$ has period $T^+$ and 
$\gamma^-_i \in \mathcal{P}_\alpha$ has period $T^-_i$ for $i= 1, \dots, r$. If the moduli space $\mathcal{M}(\gamma^+ ; \gamma^-_1, 
\dots, \gamma^-_r)$ is nonempty, then by Stokes theorem $0 \le A(u) = T^+ - \sum_{i=1}^r T^-_i$ with equality if and only if $r=1$ and 
$\gamma^+ = \gamma^-_1$.

\subsection*{\textup{2.3.} Moduli spaces of tree-like holomorphic curves} 

Let us now turn to the properties of the moduli spaces $\mathcal{M}(\gamma^+ ; \gamma^-_1, \dots, \gamma^-_r)$ of tree-like curves. In this
section, we will refer to the original articles for various facets in the study of these moduli spaces, though most of these topics are covered in the
very comprehensive lectures by 
Wendl
\cite{arxiv:1612.01009}.
The local 
structure of these moduli spaces can be studied using tools from functional
analysis and, more specifically, Fredholm theory. Our moduli space is
seen as a subset of a larger configuration space, consisting of all maps $u$ defined on a sphere with $r+1$ punctures and 
satisfying the asymptotic conditions corresponding to $(\gamma^+ ; \gamma^-_1, \dots, \gamma^-_r)$. In order to make use of Fredholm theory,
this configuration space must be equipped with a suitable Banach structure, and the standard choice is to require our maps to satisfy some 
Sobolev regularity. The minimal choice is to require one weak derivative, since
the Cauchy--Riemann equation has order $1$. For the Sobolev 
regularity to make sense for maps between manifolds, these maps must be at least continuous, and for the Sobolev embedding theorem to 
apply for our maps $u$ defined on surfaces, we must work with exponent $p > 2$. This leads to $W^{1,p}$-type Sobolev spaces. Note that if one 
prefers to work in the Hilbert case $p=2$, one can instead require two weak derivatives and use $W^{2,2}$-type Sobolev spaces instead. This 
latter choice was made in the literature corresponding to some other variants of
contact homology, but here we will keep 
to the first choice, so
that our configuration space $\mathcal{B} = \mathcal{B}(\gamma^+ ; \gamma^-_1, \dots, \gamma^-_r)$ will be a Banach manifold modeled 
on a $W^{1,p}$ space.

Let us specify the measure of the punctured sphere which is used to define the Sobolev norm. In Floer theory, one uses the standard measure 
$ds \, d\theta$ on the cylinder $\mathbb{R} \times S^1$ with coordinates $(s,\theta)$. It is therefore natural in the case of a punctured sphere 
$\mathbb{C}\mathrm{P}^1 \setminus \{ x_1, \dots, x_q \}$ to use a measure that coincides with $ds \, d\theta$ on the punctured
neighborhoods $\mathbb{R}^- \times S^1$ 
of $x_1, \dots, x_q$ with coordinates $(s, \theta) = (\log \frac{\rho}{\rho_0}, \theta)$, where $(\rho, \theta)$ are the polar coordinates used in 
Section~2.2 
and $\rho < \rho_0$ in these neighborhoods.
There is however a complication in the case of symplectizations, 
as compared to that of 
Floer theory. Despite the fact that closed Reeb orbits are 
assumed to be nondegenerate (which is a statement about the linearized return map restricted to $\xi$), the linearized return map along a closed 
orbit in the whole tangent space $T(\mathbb{R} \times M)$ of the symplectization
has an eigenvalue 1 with, as corresponding eigenspace, the plane spanned
by the Reeb vector field $R_\alpha$ and the Liouville vector field
$\frac{\partial}{\partial t}$. In other words, nondegenerate closed Reeb orbits,
seen
as Hamiltonian periodic orbits in the symplectization $\mathbb{R} \times M$, are degenerate in the sense of Floer theory. This will cause some analytic 
complications, and in order to avoid these it is necessary that 
we add exponential weights to our measure near the punctures, and use the measure
$e^{d|s|} \,ds \, d\theta$ in the corresponding neighborhoods. The corresponding Sobolev space with one weak derivative and exponent $p$ are
the referred to as a $W^{1,p,d}$ space. It can be shown that holomorphic maps $u$ with finite Hofer energy will converge with exponential speed 
to vertical cylinders over closed Reeb orbits near the punctures; 
see 
\cite{mr:1395676}
for estimates in dimension 
three and 
\cite[Appendix~A]{mr:2026549}
for the setting to generalize these to higher dimensions. More precisely, in the above $(s,\theta)$ coordinates around a puncture where 
$u_\mathbb{R} \to \pm\infty$, there exists a closed Reeb orbit $\gamma$, as well as constants $s_0 \in \mathbb{R}$ and 
$\theta_0 \in \mathbb{R}/2\pi\mathbb{Z}$ such that 
\begin{equation} \label{eq:expconv1}
	| u_\mathbb{R} (s,\theta) \pm T(s - s_0)| \le C e^{-d|s|} \tag{2.1}
\end{equation}
and
\begin{equation} \label{eq:expconv2}
	d_M(u_M(s,\theta), \gamma(\mp T(\theta-\theta_0)/2\pi)) \le C e^{-d|s|},
	\tag{2.2}
\end{equation}
using some auxiliary distance $d_M$ in $M$, for $C > 0$ sufficiently large and $d > 0$ sufficiently small. Therefore, the configuration space 
$\mathcal{B}$ defined using a $W^{1,p,d}$ space will contain the holomorphic curves we are interested in and is a suitable choice for 
constructing in it the desired moduli space. 

Note that, in 
Equations \eqref{eq:expconv1} and \eqref{eq:expconv2}, the constants $s_0$ and $\theta_0$ depend on the map $u \in \mathcal{B}$.
Because of this, the tangent space $T_u \mathcal{B}$ will not only consist of sections of $u^*T(\mathbb{R} \times M)$ that are in $W^{1,p,d}$, but also of 
sections which converge near punctures to nonvanishing vectors in the plane spanned by the Reeb field $R_\alpha$ and the Liouville field 
$\frac{\partial}{\partial t}$. Therefore, in addition to a $W^{1,p,d}$ space of sections of $u^*T(\mathbb{R} \times M)$, the tangent space $T_u \mathcal{B}$
will also contain a 
two-dimensional summand for each puncture, spanned by sections supported in a neighborhood of the puncture and taking
the constant value $R_\alpha$ or $\frac{\partial}{\partial t}$ in a smaller neighborhood of the puncture. The total dimension of these extra 
summands in $T_u \mathcal{B}$ is $2q$, where $q$ is the number of punctures. The norm of these special sections will be
determined by the asymptotic value at the corresponding puncture via an auxiliary metric on $\mathbb{R} \times M$.

In additional to these two-dimensional summands, the tangent space $T_u \mathcal{B}$ will also contain a finite dimensional summand 
corresponding to the change of position for the punctures in
$\mathbb{C}\mathrm{P}^1$. Each puncture introduces 
two degrees of freedom, but the quotient
by action of biholomorphisms of $\mathbb{C}\mathrm{P}^1$, acting transitively on triplets of points in
$\mathbb{C}\mathrm{P}^1$, reduces these degrees of freedom by $6$, so
we obtain $2(q-3)$ degrees of freedom corresponding to the position of the $q$ punctures.

For a map $u\in\mathcal{B}$, the expression $du - J \circ du \circ j$ is a $(0,1)$-form with $L^p$ regularity on the domain of $u$ 
with values in $u^*T(\mathbb{R} \times M)$. One therefore reinterprets the
Cauchy--Riemann equation as a section of a Banach bundle 
$\mathcal{E} \to \mathcal{B}$ with fibers $\mathcal{E}_u$ that are modeled on an $L^{p,d}$ space of suitable vector-valued $(0,1)$-forms 
over the configuration space $\mathcal{B}$. 
The moduli space is then realized as the zero set $s^{-1}(0)$ of such a section $s$. In the case of a finite rank vector bundle over a finite 
dimensional manifold, it suffices to
show that the section is sufficiently regular and is a submersion along its zero set in order to show that $s^{-1}(0)$ is a manifold of a given 
dimension. In the above Banach setting, the regularity condition on $s$ is replaced by the Fredholm property: the vertical differential 
$D_u : T_u \mathcal{B} \to \mathcal{E}_u$ of $s$ at some map $u$, which is a
linear, bounded operator between Banach spaces, should 
have a finite dimensional kernel, a closed image and a cokernel, 
i.e., the quotient of its target space by its image, of finite dimension as well. In that case, the Fredholm index of $D$ defined as 
$\operatorname{ind} D_u = \dim \ker D_u - \dim \operatorname{coker} D_u$ is locally constant on the space of Fredholm operators. In our case, it will
depend only on the combinatorial data decorating the moduli space under consideration, so that it does not depend on $u \in \mathcal{B}$. 
Furthermore, if the differential of $s$ is surjective along $s^{-1}(0)$, then this zero locus is a smooth manifold of finite dimension given by 
$\operatorname{ind} D_u$. This Fredholm index for a Cauchy--Riemann type operator on a Riemann surface can be computed via the 
Riemann--Roch theorem, which says that the index of a Cauchy--Riemann operator on sections of a vector bundle of complex rank $n$ 
and of first Chern class $c_1$ over a closed Riemann surface of genus $g$ is given by $n(2-2g) + 2c_1$. There are several ways to compute 
the index of $D_u$, but here is the author's favorite 
approach. First, restrict the operator $D_u$ to the $W^{1,p,d}$ part of its
domain, or, in other words, remove from its domain the summands 
of total dimension $4q -6$ corresponding to the variations of $s_0$ and
$\theta_0$ in 
Equations \eqref{eq:expconv1} and \eqref{eq:expconv2}, as well
as the positions of the punctures. Second, conjugate the restricted operator with the multiplication by a function of the form $e^{d|s|}$ near each
puncture, so that the resulting operator is defined in a $W^{1,p}$ space and takes its values in an $L^p$ space. While the original operator has the 
form $D_u \zeta = \partial_s \zeta + J_0 \partial_\theta \zeta + A_j(s, \theta) \zeta$ near puncture $j$, in a
suitable trivialization of $u^*T(\mathbb{R} \times M)$,
the conjugated operator will have a similar form with the matrix $A_j$ replaced with $\tilde{A}_j = A_j - \varepsilon_j d I$ where $\varepsilon_j$ is 
the sign of the puncture $j$, with $j = 1, \dots, q$. 
Near the punctures, such operators are exactly of those encountered as the linearization of the Floer equation in symplectic geometry. One can
show 
\cite[Theorem~3.2.12]{cm:schwarz.phdthesis.1995}
that the index of these operators is additive under the gluing operation at one or several punctures. 
The trick is then to glue 
two identical operators along corresponding pairs of punctures, using a suitable clutching function for the vector bundle at each puncture so that the 
expressions of the operators match in the regions that are to be identified when gluing the Riemann surfaces. The glued Riemann surface 
has genus $q-1$ and the vector bundle over it has its first Chern class given by $\sum_{j=1}^q \varepsilon_j \mu_{CZ}(\tilde{A}_j)$, where 
$\varepsilon_j$ is the sign of the puncture $j$. By the Riemann--Roch theorem, this glued operator has index 
\[
	n(2-2(q-1)) + 2 \sum_{j=1}^q \varepsilon_j \mu_{CZ}(\tilde{A}_j),
\]
and by the additivity of the index, this is twice the index of the conjugated 
operator. A simple calculation shows that $\mu_{CZ}(\tilde{A}_j) = \mu_{CZ}(A_j) - \varepsilon_j$, so that the index of the restricted operator
can be expressed as 
\[
	n(2-q) - q + \sum_{j=1}^q \varepsilon_j \mu_{CZ}(A_j).
\]
Adding the removed degrees of freedom, the index of $D_u$
is given by 
\[
n(2-q) + 3q - 6 + \sum_{j=1}^q \varepsilon_j \mu_{CZ}(A_j).
\]
Since the $q$ punctures correspond to one positive puncture and 
$r$ negative punctures, and the indices $\mu_{CZ}(A_j)$ for $j = 1 , \dots, q$  are given by $\mu_{CZ}(\gamma^+)$ and $\mu_{CZ}(\gamma^-_i)$ 
for $i = 1, \dots, r$, we finally obtain
\begin{equation} \label{eq:index}
\operatorname{ind} D_u = (n-3)(1-r) + \mu_{CZ}(\gamma^+) - \sum_{i=1}^r
	\mu_{CZ}(\gamma^-_i).\tag{2.3}
\end{equation}
This formula gives the dimension of the moduli space $\widetilde{\mathcal{M}}(\gamma^+ ; \gamma^-_1, \dots, \gamma^-_r)$, provided 
one can
ensure that the operator $D_u$ is surjective for all $[u]$ in this moduli space.

Under this surjectivity assumption, one can construct a gluing map with a sufficiently large parameter $R > 0$ which is defined over 
any compact subset of $\widetilde{\mathcal{M}}(\Gamma_1) \times \widetilde{\mathcal{M}}(\Gamma_2)$, where the last negative asymptote 
of $\Gamma_1$ coincides with the positive asymptote of $\Gamma_2$. This map takes its values in $\widetilde{\mathcal{M}}(\Gamma)$, 
where $\Gamma$ has the same positive asymptote as $\Gamma_1$, and its negative asymptotes are the negative asymptotes of $\Gamma_1$ 
(except the last one) and the negative asymptotes of $\Gamma_2$. This map takes a pair of holomorphic curves $(u_1, u_2)$ to a holomorphic 
curve $u_R$ that nearly coincides with $(u_{1, \mathbb{R}} + R, u_{1,M})$ on $\mathbb{R}^+ \times M$ and with
$(u_{2,\mathbb{R}} -R, u_{2,M})$ on $\mathbb{R}^- \times M$. 
When the positive asymptote of $\Gamma_2$  is a closed Reeb orbit covering $m$ times a simple closed orbit, there are $m$ different ways 
of gluing $u_1$ and $u_2$, differing by a shift by any multiple $2\pi/m$ of the coordinate $\theta$ near the positive puncture of $u_2$ before
identification of this coordinate with an analogous coordinate near the last negative puncture of $u_1$. This combinatorial choice must be 
specified as part of the gluing data for the moduli spaces $\widetilde{\mathcal{M}}(\Gamma_1)$ and $\widetilde{\mathcal{M}}(\Gamma_2)$.
The construction of this gluing map follows the same strategy as other uses of holomorphic curves in contact and symplectic geometry, but relies 
on the specific Banach structures for contact homology that are described above. One first constructs a preglued map $\tilde{u}_R$ which is 
constructed from the maps $u_1$ and $u_2$ using cutoff functions near their common asymptote.
Then 
one applies an infinite dimensional version of the implicit function theorem to find a holomorphic curve $u_R$ near $\tilde{u}_R$ in the configuration
space $\mathcal{B}(\Gamma)$. The application of this theorem requires three estimates. First, the $L^{p,d}$ norm of 
$d\tilde{u}_R - J \circ d\tilde{u}_R \circ j$ must tend to $0$ as $R \to \infty$. Second, the operator $D_{\tilde{u}_R}$ must have a right inverse which
is uniformly bounded in $R$, with respect to the Banach norms described above. Third, the second order term in the Taylor expansion of the 
section $s : \mathcal{B} \to \mathcal{E}$ around $\tilde{u}_R$ must be bounded uniformly in $R$ with respect to the same norms. The second 
estimate is typically the most delicate one to establish. To this end, one first constructs an approximate right inverse $Q_R$ for $D_{\tilde{u}_R}$ 
using bounded right inverses for $D_{u_1}$ and $D_{u_2}$, this is where the surjectivity assumption mentioned above is essential. More precisely,
one has to prove that the norm of the operator $D_{\tilde{u}_R} \circ Q_R - I$ tends to $0$ as $R \to \infty$. Then, standard arguments provide
an actual right inverse with the desired properties for $R$ sufficiently large.
Once this gluing map is defined, it can be used to construct the compactification of the moduli spaces of holomorphic curves as a smooth manifold
with boundaries and corners. Indeed, it turns out that holomorphic buildings with $\ell$ levels constitute a codimension $\ell - 1$ submanifold in the
boundary of the compact moduli space of holomorphic buildings.

In order to define orientations having suitable properties with respect to the above gluing map on our moduli spaces, it is necessary 
to work in a more abstract setting. Let $\mathcal{O} = \mathcal{O}(\gamma^+ ; \gamma^-_1, \dots, \gamma^-_r)$ be the space of 
differential operators which are of the same form as the above operators $D_u$ with $u \in \mathcal{B}$. In particular, all 
$D \in \mathcal{O}$ are Fredholm operators with
the same index as above. Although the dimensions of the vector spaces $\ker D$
and $\operatorname{coker} D$ can vary as $D$ varies
in $\mathcal{O}$, the real lines $\Lambda^{\max} \ker D \otimes \Lambda^{\max} (\operatorname{coker} D)^*$ naturally fit together
to form a real line bundle $\mathcal{L} \to \mathcal{O}$, called 
a \textit{determinant line bundle}; this is true for general spaces of Fredholm 
operators; see 
\cite[Appendix]{mr:1200162}.
An orientation of $\mathcal{L}$ is equivalent to the data of a nonvanishing section of $\mathcal{L}$.
Such a section exists if and only if $\mathcal{L}$ is a trivial vector bundle. Note that $\mathcal{O}$ is homotopy equivalent to a product
of $q=r+1$ circles, because of the possible values of $\theta_0 \in
\mathbb{R}/2\pi\mathbb{Z}$ in 
Equation \eqref{eq:expconv2}. It turns out that $\mathcal{L}$ is
trivial along a circle factor corresponding to a puncture where $D$ has the asymptotic behavior for the asymptote $\gamma$ if and only
if $\gamma$ is a good orbit 
\cite[Theorem~3]{mr:2092725}.
This is the reason for the
distinction between good and bad orbits in Section~2.1. 
If we restrict ourselves to good orbits as asymptotes of holomorphic curves, then these considerations guarantee that all moduli spaces are
orientable. One can then construct 
\cite{mr:2092725}
a coherent set of orientations for the determinant line bundles corresponding to all operator spaces 
$\mathcal{O}(\Gamma)$. Here the word ``coherent'' means that these orientations are preserved by the gluing maps 
$\mathcal{L}(\Gamma_1) \otimes \mathcal{L}(\Gamma_2) \to \mathcal{L}(\Gamma)$ induced by the gluing of operators in a similar way as the
gluing map for holomorphic curves discussed above. These coherent orientations are uniquely determined by choices of orientations for the
determinant line bundles with a single (positive) asymptote at all good closed Reeb orbits. One subtle point is that the coherent orientation 
on $\mathcal{L}(\gamma^+; \gamma^-_1, \dots, \gamma^-_r)$ also depends on the ordering of the negative asymptotes. This can be understood
by the following considerations: gluing the above coherent orientation successively with the coherent orientations on 
$\mathcal{L}(\gamma^-_i)$ for $i= r, \dots, 1$, one obtains the coherent orientation on $\mathcal{L}(\gamma^+)$. But the coherent orientation
on $\mathcal{L}(\gamma^-_r) \otimes \dots \otimes \mathcal{L}(\gamma^-_1)$ depends on the ordering of the oriented bases for the kernel 
and cokernel of operators in $\mathcal{O}(\gamma^-_i)$ for $i = r, \dots, 1$. Permuting two consecutive asymptotes $\gamma^-_i$ and 
$\gamma^-_{i+1}$ will change the product orientation exactly when both permuted bases consist of an odd number of vectors, or in other 
words when the Fredholm indices \eqref{eq:index} of the corresponding operators $\mu_{CZ}(\gamma^-_i)+n-3$ and 
$\mu_{CZ}(\gamma^-_{i+1})+n-3$ are both odd. Pulling back this system of coherent orientations on the determinant line bundles by the
map $\widetilde{\mathcal{M}}(\Gamma) \to \mathcal{O}(\Gamma) : u \mapsto D_u$
gives a coherent system of orientations on the moduli spaces
$\widetilde{\mathcal{M}}(\Gamma)$. In the case of a surjective operator $D_u$, its determinant line is indeed given by 
$\Lambda^{\rm max}(\ker D_u)$, where $\ker D_u$ is naturally identified with the tangent space $T_u\widetilde{\mathcal{M}}(\Gamma)$.
Note that this system of coherent orientations is defined for all moduli spaces, regardless of their dimensions. In addition to this, one can 
define canonical 
orientations 
for all moduli spaces $\widetilde{\mathcal{M}}(\Gamma)$ of dimension 
one. These are defined by pushing forward 
the natural orientation of the real line by the $\mathbb{R}$-action on
$\widetilde{\mathcal{M}}(\Gamma)$. Then all one-dimensional moduli spaces are
equipped with two orientations: the coherent one and the canonical one. Each connected component of $\widetilde{\mathcal{M}}(\Gamma)$, or
equivalently each element $[u]$ of $\mathcal{M}(\Gamma)$, can therefore be equipped with a sign: $+1$ if both orientations agree, and $-1$
otherwise.
 
The above construction of the moduli spaces $\mathcal{M}(\Gamma)$ as oriented finite dimensional manifolds relies on the assumption that 
the operator $D_u$ is surjective for every $[u] \in \mathcal{M}(\Gamma)$. In other theories involving holomorphic curves in contact and symplectic
topology, this property is achieved using so-called classical transversality techniques 
\cite{mr:2045629},
which consist in choosing the almost complex 
structure $J$ generically. However, these methods fail in the presence of
multiply covered holomorphic curves, or, 
in other words, 
holomorphic 
curves that factor though a ramified covering of Riemann 
surfaces. Such covering can unavoidably occur, for example, in the case of a moduli space $\mathcal{M}(\gamma)$ of holomorphic planes with
a single positive puncture asymptotic to a closed Reeb orbit $\gamma$ which is
the iterate of another closed Reeb orbit. Several other methods have been 
developed to overcome these very delicate technical difficulties. Most of
them consist in perturbing the right-hand side of the Cauchy--Riemann 
equation. When multiply covered holomorphic curves arise, these come with a nontrivial group of automorphisms, and it is typically impossible to
achieve transversality of the section $s : \mathcal{B} \to \mathcal{E}$ with the zero section at such curves using an equivariant right-hand side. For
this reason, it is necessary to used so-called multivalued perturbations, or, in
other words, a finite set of right-hand sides for the Cauchy--Riemann 
equation which is globally preserved by the automorphism group. When this set consists of more than one element, it is important not to overcount 
holomorphic curves, so that the perturbations, as well as the corresponding elements in the perturbed moduli space, have to carry fractional weights 
so that the sum of all relevant weights remains equal to $1$. This is why the
count $n(\Gamma)$ of elements in a perturbed 
zero-dimensional 
moduli space $\mathcal{M}(\Gamma)$ is a rational number, as each element comes both with a sign coming from the orientations and with a 
fractional weight coming from the necessary perturbation to achieve transversality. When no such perturbation is required, $n(\Gamma) \in 
\mathbb{Q}$ is the sum over all elements $[u] \in \mathcal{M}(\Gamma)$ of the sign of $[u]$ divided by the order of the automorphism 
group of $u$.

Beyond this very brief description of the perturbation scheme for our moduli
spaces, the complete constructions are extremely long and difficult to
describe 
in great detail.
Over the years, many different research teams have developed their own approach to this very complicated problem. 
The theory of polyfolds introduced by Hofer, 
Wysocki
and 
Zehnder
\cite{mr:4298268}
aims to frontally attack 
all analytical problems in order to produce moduli spaces that are as close as possible to honest manifolds with boundaries and corners.
Kuranishi structures were used by 
Fukaya
and 
Ono
\cite{mr:1688434}
in the context of Gromov--Witten invariants and these may be used for 
contact homology in order to obtain suitable moduli spaces. In the same spirit, 
Bao
and  
Honda
\cite{mr:4539062}
used semiglobal Kuranishi 
charts to provide a definition of the contact homology algebra. 
Hutchings
and 
Nelson
\cite{mr:3601881},
\cite{mr:4461852}
used automatic transversality 
techniques to define contact homology in dimension three. 
Pardon
\cite{mr:3493097},
\cite{mr:3981989}
developed an algebraic, rather than geometric, approach
via virtual fundamental cycles, in order to define suitable counts of elements in low-dimensional moduli spaces, and applied this technique
to define contact homology in full generality.

\section*{3. Contact homology algebra}

In this section, we introduce the algebraic framework that mirrors the geometric properties of the moduli spaces of holomorphic curves
and leads to the definition of homological invariants for contact manifolds, namely contact homology and its variants.

\subsection*{\textup{3.1.} Differential graded algebra}

Since the moduli spaces of holomorphic curves from Section~2.2 
are decorated with closed Reeb orbits, it is necessary to
make these orbits part of the algebraic formalism that will incorporate the count of elements in these moduli spaces. To each closed orbit
$\gamma \in \mathcal{P}^{\rm g}_\alpha$, we associate a formal generator $q_\gamma$, in order to distinguish geometric objects from algebraic 
ones. The generator $q_\gamma$ will be given a grading obtained by looking at the dimension formula for moduli spaces of holomorphic 
curves. We set $|q_\gamma| = \mu_{CZ}(\gamma) + n - 3$, so that in view of
Equation \eqref{eq:index} we have
\begin{equation} \label{eq:dim}
\dim \mathcal{M}(\gamma^+; \gamma^-_1, \dots, \gamma^-_r) = |q_{\gamma^+}| -
	\sum_{i=1}^r |q_{\gamma^-_i}| -1.\tag{3.1}
\end{equation}

Let $\mathcal{A}$ be the supercommutative, graded, unital algebra generated by $q_\gamma$ for all $\gamma \in \mathcal{P}^{\rm g}_\alpha$
over some ring $\mathcal{R}$. We are working here with an algebra, and not a module as in Floer theory, because our holomorphic curves are
allowed to have several negative asymptotes (i.e., when $r > 1$), and in this case we will be using the formal expression 
$q_{\gamma^-_1} \dots q_{\gamma^-_r}$. Moreover, we need a unital algebra so that the unit $1 \in \mathcal{A}$ can be used when there are
no negative asymptotes (i.e., when $r = 0$). In the above definition,
``supercommutative'' means that for any pair of elements $a, b \in \mathcal{A}$ 
with pure grading (so that $|a|$ and $|b|$ are well defined), we have $a b = (-1)^{|a|  |b|} b a$. This rule is introduced in order to mirror the 
behavior of the coherent orientations of the moduli spaces $\mathcal{M}(\gamma^+; \gamma^-_1, \dots, \gamma^-_r)$ under the exchange 
of two orbits $\gamma^-_i$ and $\gamma^-_j$.

Let us now discuss the possible choices for the ring $\mathcal{R}$. Making the choice $\mathcal{R} = \mathbb{Z}/2\mathbb{Z}$ would 
eliminate all sign considerations, but such coefficients are incompatible with the signed and weighted counts of elements in moduli spaces,
which are rational numbers, as explained in Section~2.3. 
Therefore, the simplest choice is to take
$\mathcal{R} = \mathbb{Q}$, under our assumption that $c_1(\xi)=0$. When this assumption is not satisfied, a more elaborate choice of ring
is required if we wish to retain a grading in our algebraic formalism. One then usually chooses $\mathcal{R} = \mathbb{Q}[H_2(M)]$. 
A general element of this group ring has the form $\sum_{A \in H_2(M)} q_A e^A$, where $q_A \in \mathbb{Q}$ and $e^A$ is a formal
generator with grading $|e^A| = -2 \langle c_1(\xi), A \rangle$. Once again, this grading is chosen in view of the dimension formula
for the moduli spaces of holomorphic curves. Note that the grading of $\mathcal{R}$ is even, so that it does not interfere with the 
supercommutativity property of $\mathcal{A}$. Of course, one can make this choice $\mathcal{R} = \mathbb{Q}[H_2(M)]$ even in the case
where $c_1(\xi)=0$, and, more 
generally, one can choose $\mathcal{R}$ to be any graded algebra over the graded ring 
$\mathbb{Q}[H_2(M)]$. Such special choices can be useful in certain
situations.

Let us now turn to the definition of the differential $\partial : \mathcal{A} \to \mathcal{A}$. We require $\partial$ to be a linear superderivation 
of degree $-1$ of $\mathcal{A}$, or, 
in other words, $\partial (\lambda a+ \mu b) = \lambda \partial a + \mu \partial b$ for all $a, b \in \mathcal{A}$, 
$\lambda, \mu \in \mathcal{R}$, and $|\partial a | = |a|-1$ as well as $\partial(ab) = (\partial a)b + (-1)^{|a|} a \partial b$ for all $a, b \in \mathcal{A}$ 
of pure grading. The  differential is therefore characterized by its value on the generators of $\mathcal{A}$, and we define
\begin{equation} \label{eq:d}
\partial q_{\gamma^+} = m_{\gamma^+} \sum_{r \ge 0} 
\sum_{\substack{\{\gamma^-_1, \dots, \gamma^-_r\}  \subset \mathcal{P}^g_\alpha \\ |q_{\gamma^-_1}| + \dots + |q_{\gamma^-_r}| 
= |q_{\gamma^+}|-1}}
n(\gamma^+;\gamma^-_1, \dots, \gamma^-_r) \ q_{\gamma^-_1} \dots
	q_{\gamma^-_r},\tag{3.2}
\end{equation}
where $m_{\gamma^+}$ is the multiplicity of the orbit $\gamma^+$ and $n(\gamma^+;\gamma^-_1, \dots, \gamma^-_r) \in \mathbb{Q}$ is the signed 
and weighted count of the elements $[u]$ in the moduli space
$\mathcal{M}(\gamma^+;\gamma^-_1, \dots, \gamma^-_r)$, 
as described in 
Section~2.3. 
Note that this sum makes sense because our definition of the grading 
in $\mathcal{A}$ in terms of the dimension formula for the moduli spaces guarantees that the moduli spaces involved in this sum 
are zero-dimensional, and hence finite by compactness.
Also note that, in view of the supercommutativity rule in $\mathcal{A}$ mimicking the behavior of coherent orientations for moduli spaces, the above
expression is well defined independently of the ordering of the generators $q_{\gamma^-_1}, \dots, q_{\gamma^-_r}$. 

It then follows from the same philosophy as in Morse or Floer theory, using the properties of the moduli spaces of holomorphic curves, that 
we have the identity $\partial \circ \partial = 0$. In other words, $(\mathcal{A}, \partial)$ is a differential graded algebra (DGA). More precisely, 
this is proved by 
considering the boundary of one-dimensional moduli spaces $\mathcal{M}(\gamma^+;\gamma^-_1, \dots, \gamma^-_r)$. The signed weighted 
count of the elements in this boundary then corresponds to the coefficient of $q_{\gamma^-_1} \dots q_{\gamma^-_r}$ in the expression
$\partial \circ \partial q_{\gamma^+}$. But the sum of the weights associated to elements in the boundary of a 
connected component of the above one-dimensional moduli space vanishes 
by definition of the weights and of the signs. This implies
that $\partial \circ \partial q_{\gamma^+} = 0$ for any generator $q_{\gamma^+}$.

The homology $\ker \partial / \operatorname{im} \partial$ of the DGA $(\mathcal{A}, \partial)$ is a graded 
algebra; it is called the \textit{contact 
homology} of $(M, \xi)$ and is denoted by $CH(M,\xi)$. This terminology and this notation are justified by the fact that, unlike the 
DGA $(\mathcal{A}, \partial)$, contact homology does not depend on the various choices made during this construction, in particular 
the contact form $\alpha$ for $\xi$ and the compatible complex structure $J$ on $(\xi, d\alpha)$.

The proof of the invariance property of contact homology is similar in philosophy to the proof of $\partial \circ \partial = 0$. The main 
difference lies in the fact that we are not considering holomorphic curves in a symplectization anymore, but in a symplectic cobordism
of the form $(\mathbb{R} \times M, d(e^t \alpha_t))$, with $\alpha_t$ interpolating very slowly between different choices of contact forms: 
$\alpha_t = \alpha_+$ for $t$ sufficiently 
large, and $\alpha_t = \alpha_-$ for $t$ sufficiently small. Furthermore, the fact that $\alpha_t$
varies slowly with $t$ guarantees that the $2$-form $d(e^t \alpha_t)$ is nondegenerate, hence symplectic. This symplectic cobordism is 
also equipped with a compatible almost complex structure that depends on $t$ and interpolates between a compatible almost complex 
structure $J_+$ for the symplectization of $\alpha_+$ and a compatible almost complex structure $J_-$ for the symplectization of $\alpha_-$.
Different discrete choices can also be made at the ends of this cobordism in order to define the corresponding DGAs 
$(\mathcal{A}_+, \partial_+)$ and $(\mathcal{A}_-, \partial_-)$ via their respective symplectizations.

The properties of the moduli spaces of holomorphic curves in such a symplectic cobordism are very similar to those in the case of a 
symplectization, except for the fact that they are not equipped with a free action of $\mathbb{R}$ by translation along the first coordinate of the
cobordism. Therefore, the dimension formula for these moduli spaces does not
contain the final $-1$ term as in Equation \eqref{eq:dim}, and the 
count of elements in zero-dimensional moduli spaces will lead to the definition of a map of degree $0$ instead of $-1$. This map 
$\Phi : (\mathcal{A}_+, \partial_+) \to (\mathcal{A}_-, \partial_-)$ is defined
to be a morphism of unital algebras, in other words, $\Phi(1) = 1$,
$\Phi(\lambda a + \mu b) = \lambda \Phi(a) + \mu \Phi(b)$ and $\Phi(ab) = \Phi(a) \Phi(b)$ for all $a, b \in \mathcal{A}_+$ and 
$\lambda, \mu \in \mathcal{R}$. This map is then also characterized by its value on generators $q_{\gamma^+}$ of $\mathcal{A}_+$, 
and a formula similar to Equation \eqref{eq:d} is used to define $\Phi(q_{\gamma^+})$.

The boundary of one-dimensional moduli spaces of holomorphic curves in the above symplectic cobordism consists of holomorphic buildings 
with two levels, one in the cobordism and one in the symplectization of $\alpha_+$ or of $\alpha_-$. Therefore, the vanishing weighted count
of the elements in this boundary leads to the identity $\Phi \circ \partial_+ = \partial_- \circ \Phi$. In other words, $\Phi$ is a unital 
DGA map and induces a map $\overline{\Phi}$ between the homologies of $(\mathcal{A}_+, \partial_+)$ and $(\mathcal{A}_-, \partial_-)$. It remains
to show that this induced map is actually an isomorphism.

Exchanging the roles of $(\mathcal{A}_+, \partial_+)$ and $(\mathcal{A}_-, \partial_-)$, and reversing the sign of $t$ in the interpolations 
$\alpha_t$ and $J_t$, one obtains similarly a unital DGA map $\Psi : (\mathcal{A}_-, \partial_-) \to (\mathcal{A}_+, \partial_+)$ inducing
a map $\overline{\Psi}$ between the homologies of $(\mathcal{A}_-, \partial_-)$ and $(\mathcal{A}_+, \partial_+)$. A gluing argument then
shows that the composition $\overline{\Psi} \circ \overline{\Phi}$ is induced by the concatenation of both interpolations, which is isotopic to
the constant interpolation between $(\alpha_+, J_+)$ and itself. Note that this constant interpolation induces the identity maps on 
$(\mathcal{A}_+, \partial_+)$ and on its homology, because the moduli spaces
that are zero-dimensional without taking the quotient
by the $\mathbb{R}$-action on the first coordinate of the symplectization consist purely of vertical cylinders over closed Reeb orbits. If one could
show that the map induced in homology by an interpolation as above is invariant under isotopy of such interpolations, then we could deduce
that $\overline{\Psi} \circ \overline{\Phi}$ is the identity. Similarly, $\overline{\Phi} \circ \overline{\Psi}$ would also be the identity, so that
$\overline{\Phi}$ would be an isomorphism as desired.

Given an isotopy of interpolations between interpolations inducing unital DGA maps 
\[
	\Phi_0, \Phi_1 : (\mathcal{A}_+, \partial_+) \to (\mathcal{A}_-, \partial_-),
\]
we can construct a one-parameter family of symplectic cobordisms
of the form $(\mathbb{R} \times M, d(e^t \alpha_{t,\sigma}))$ 
equipped with compatible almost complex structures $J_\sigma$, for $\sigma \in [0,1]$.
We consider the moduli space $\mathcal{M}^\sigma(\gamma^+;\gamma^-_1, \dots, \gamma^-_r)$ of holomorphic curves in this family of cobordisms,
consisting of all pairs $(\sigma, [u_\sigma])$ where $\sigma \in [0,1]$ and $u_\sigma$ is a $J_\sigma$-holomorphic map in 
$(\mathbb{R} \times M, d(e^t \alpha_{t,\sigma}))$. By the compactness theorem
for holomorphic curves in symplectic cobordisms, the zero-dimensional
moduli spaces $\mathcal{M}^\sigma$ are finite, 
so that only finitely many values $\sigma_1, \dots, \sigma_N \in [0,1]$ appear as the 
first coordinate of an element in such a moduli space. If $N=0$, then the
boundary elements of one-dimensional moduli spaces $\mathcal{M}^\sigma$ 
must be of the form $(0, [u_0])$ or $(1,[u_1])$, so that $\Phi_0 = \Phi_1$. We can therefore concentrate on small intervals around the special values
$\sigma_i$ for $i= 1, \dots, N$ or, in other words, assume that $N=1$. A $(\Phi_0, \Phi_1)$-derivation $K : \mathcal{A}_+  \to \mathcal{A}_-$
is a map satisfying $K(1) = 0$, $K(\lambda a + \mu b) = \lambda K(a) + \mu K(b)$ for all $a, b \in \mathcal{A}$ and $\lambda, \mu \in \mathcal{R}$, 
and
\begin{align*}
K(q_1 \dots q_k) &= \sum_{g \in S_k} \sum_{j=1}^k (-1)^{|q_{g(1)} \dots q_{g(j-1)}|} 
\Phi_0(q_{g(1)} \dots q_{g(j-1)}) K(q_{g(j)}) \Phi_1(q_{g(j+1)} \dots q_{g(k)}),
\end{align*}
for all generators $q_1, \dots, q_k$ of $\mathcal{A}_+$, where $S_k$ denotes the permutation group of $\{ 1, \dots, k \}$.
We define $K(q_{\gamma^+})$ by a formula similar to Equation \eqref{eq:d} using
the signed and weighted count of elements in the zero-dimensional moduli spaces 
$\mathcal{M}^\sigma(\gamma^+;\gamma^-_1, \dots, \gamma^-_r)$. Since the parameter $\sigma \in [0,1]$ provides an additional degree 
of freedom in comparison with the previous situation, the dimension formula for these moduli spaces has a final term $+1$ instead of $-1$ 
in Equation \eqref{eq:dim}, so that $|q_{\gamma^+}| = |q_{\gamma^-_1} \dots q_{\gamma^-_r}| - 1$ and $K$ is of degree $+1$ as announced.

We then study the compactification of the one-dimensional moduli space $\mathcal{M}^\sigma(\gamma^+;\gamma^-_1, \dots, \gamma^-_r)$
using a suitable perturbation scheme. When a sequence of holomorphic curves degenerates into a holomorphic building with $2$ levels, which 
necessarily occurs when $\sigma \to \sigma_1$, it is necessary to arrange the perturbation scheme so that each connected component in the 
level corresponding to the family of symplectic cobordism has its own value of $\sigma$ very close to $\sigma_1$. Only the component counted by
$K$ has $\sigma = \sigma_1$, while the other components have either $\sigma < \sigma_1$, so that it is counted by $\Phi_0$, or $\sigma > \sigma_1$, 
so that it is counted by $\Phi_1$. With a perturbation scheme that averages over all the possible inequalities, the count of elements in the boundary of 
the one-dimensional moduli space $\mathcal{M}^\sigma(\gamma^+;\gamma^-_1, \dots, \gamma^-_r)$ leads to the identity 
$\Phi_1 - \Phi_0 = K \circ \partial_+ + \partial_- \circ K$. It follows that the induced maps on homology coincide: $\overline{\Phi}_0 = \overline{\Phi}_1$.

\subsection*{\textup{3.2.} Cylindrical contact homology}

Since the differential graded algebra in the previous section is quite large and its homology can be difficult to compute, it can be convenient to 
insist on working with holomorphic cylinders only in order to obtain an invariant which is simpler to compute. Of course, this will be possible only
under suitable conditions. 

If our holomorphic curves have a single negative puncture, the output of the differential will be a linear combination of generators corresponding 
to closed Reeb orbits, so that it is not necessary to use the algebra $\mathcal{A}$ anymore. Using the same type of ring $\mathcal{R}$ as in the
previous section, let $C^{\rm cyl}$ be the graded $\mathcal{R}$-module freely generated by $q_\gamma$ for all $\gamma \in 
\mathcal{P}^{\rm g}_\alpha$.

We define the cylindrical differential $\partial^{\rm cyl} : C^{\rm cyl} \to C^{\rm cyl}$ as the linear map of degree $-1$ characterized by
\[
\partial^{\rm cyl} q_{\gamma^+} = m_{\gamma^+} \sum_{\substack{\gamma^- \in \mathcal{P}^g_\alpha \\ |q_{\gamma^-}| = |q_{\gamma^+}|-1}}
n(\gamma^+;\gamma^-) \ q_{\gamma^-}.
\]
In other words, in the definition \eqref{eq:d} of the contact homology differential, we restrict ourselves to the term $r=1$ corresponding 
to holomorphic cylinders.

The arguments from the previous section to show that $\partial^{\rm cyl} \circ \partial^{\rm cyl} = 0$ can be adapted to holomorphic cylinders 
provided a 
one-parameter family of holomorphic cylinders can only degenerate
into a level 
two holomorphic building with a cylinder in each level.
However, there can a priori be another possible configuration for such a level
two holomorphic building: a pair of pants with one positive puncture 
and two negative punctures in the upper level, and a plane with one positive puncture together with a vertical cylinder in the lower level. In order
to prevent the existence of such a configuration, it suffices to forbid the existence of a closed Reeb orbit that can play the role of a positive 
asymptote for a rigid holomorphic plane. In view of the dimension formula \eqref{eq:dim}, we need to assume that there is no contractible 
$\gamma \in \mathcal{P}^{\rm g}_\alpha$ such that $\mu_{CZ}(\gamma) + n - 3 = 1$. Under this condition, it is possible to define cylindrical
contact homology $CH^{\rm cyl}(M, \xi)$ as the homology of the chain complex $(C^{\rm cyl}, \partial^{\rm cyl})$.

Similarly, the arguments from the previous section to show invariance of contact homology can be adapted to the cylindrical case if rigid 
holomorphic planes cannot exist in symplectic cobordisms or in one-parameter
families of symplectic cobordisms. 
Due to the dimension shifts
of the corresponding moduli spaces as in the previous section, we need to assume that there is no contractible 
$\gamma \in \mathcal{P}^{\rm g}_\alpha$ such that $\mu_{CZ}(\gamma) + n - 3$ is equal to $0$ or to $-1$. 

We can summarize the above discussion as follows: if that there is no contractible 
$\gamma \in \mathcal{P}^{\rm g}_\alpha$ such that $\mu_{CZ}(\gamma) + n - 3$ is equal to $+1$, to $0$ or to $-1$, then cylindrical
contact homology $CH^{\rm cyl}(M, \xi)$ is well defined and is an invariant of the contact manifold $(M, \xi)$. Note that more restrictive 
conditions were sometimes used in the literature to guarantee the existence of this contact invariant, when it was still an open problem to
achieve transversality for the moduli spaces of holomorphic planes. A typical condition was to require that the contact manifold $(M, \xi)$
be hypertight, or, in other words, that it admit a contact form without contractible closed Reeb orbits.

Note that, as soon as the chain complex $(C^{\rm cyl}, \partial^{\rm cyl})$ is well defined, it splits as the direct sum of chain complexes
generated by good closed Reeb orbits 
into a given free homotopy class of loops, because the differential clearly preserves those free 
homotopy classes. It can be very convenient to work with only one summand of this direct sum in order to further limit the required 
calculations to obtain a contact invariant. 

In the particular case of hypertight contact structures, if we restrict
ourselves to a free homotopy class which is primitive, or, 
in other words, not
the iterate of another free homotopy class, then all closed Reeb orbits under consideration have multiplicity one. 
This implies that the holomorphic curves that are counted by the differential of
this chain complex are not multiply covered, that is, that they
do not factor through a ramified covering of Riemann surfaces. In this case, it is actually possible to achieve transversality for the corresponding 
moduli spaces by choosing the almost complex structure $J$ generically, as shown by 
Dragnev
\cite{mr:2038115};
see also the appendix of 
\cite{mr:2200047}
for an alternative and briefer argument.

\subsection*{\textup{3.3.} Augmentations and linearization}

Although cylindrical contact homology can be a convenient alternative to the full version of contact homology, the conditions guaranteeing
that it is well defined and invariant are sometimes too restrictive and are not always natural with respect to the geometric context. In order
to overcome these limitations, one can use another method to obtain from the contact homology DGA a chain complex which is simply a 
module over the ring $\mathcal{R}$, and which is based on the notion of augmentation. Given a DGA $(\mathcal{A}, \partial)$ over a ring
$\mathcal{R}$, note that one can think of $\mathcal{R}$ as a DGA over itself, and equipped with the zero differential; it is also the
contact homology DGA of the empty contact manifold. An augmentation 
of $(\mathcal{A}, \partial)$ is defined as a unital DGA map $\varepsilon : (\mathcal{A}, \partial) \to (\mathcal{R},0)$ of degree $0$,
or in 
other words 
\[
\varepsilon(1) = 1, \quad \varepsilon(\lambda a + \mu b) = \lambda
\varepsilon(a) + \mu \varepsilon(b), \quad \varepsilon(ab) = \varepsilon(a)
\varepsilon(b),
\]
for all $a, b \in \mathcal{A}$, $\lambda, \mu \in \mathcal{R}$, and 
$\varepsilon \circ \partial = 0$.

Let us 
consider two important situations in which such augmentations arise naturally. First, if there is no contractible 
$\gamma \in \mathcal{P}^{\rm g}_\alpha$ such that $\mu_{CZ}(\gamma) + n - 3$ is equal to $+1$, to $0$ or to $-1$,
so that cylindrical contact homology is well defined, then the map $\varepsilon_{\rm cyl} : \mathcal{A} \to \mathcal{R}$
which sends $1$ to $1$ and any generator $q_\gamma$ of $\mathcal{A}$ to $0$ is an augmentation. Indeed, if $|q_{\gamma^+}|=1$, 
the constant term in $\partial q_{\gamma^+}$, corresponding to $r=0$, must vanish since the orbit $\gamma^+$ is not contractible
and can therefore not bound a holomorphic plane.

Second, if $(W, \omega)$ is a strong symplectic filling of $(M, \xi)$, or, in
other words, if a collar neighborhood of $\partial W$ is 
symplectomorphic to a portion of the symplectization of $(M, \xi)$, then one can attach the upper part of the symplectization 
of $(M, \xi)$ to $(W, \omega)$ in order to make the Liouville vector field complete. One can also extend the compatible almost 
complex structure $J$ on this part of the symplectization to the whole completed manifold. We also require that $(W, \omega)$ 
is symplectically aspherical, or, in other words, that $\omega$ evaluates trivially on $\pi_2(W)$, 
so that no bubbling of holomorphic 
spheres can occur in $(W, \omega)$. For any 
$\gamma \in \mathcal{P}_\alpha^{\rm g}$, let us denote by $\mathcal{M}_W(\gamma)$ the moduli space of $J$-holomorphic planes
in this completed manifold with one positive puncture asymptotic to the orbit $\gamma$. The dimension of this moduli space is
$|q_\gamma|$. Let $\varepsilon_W : \mathcal{A} \to \mathcal{R}$ be the unital algebra morphism of degree $0$ characterized by
$\varepsilon_W(\gamma) = n_W(\gamma)$, where $n_W(\gamma)$ is the signed and weighted count of the elements 
$[u] \in \mathcal{M}_W(\gamma)$, as described in Section~2.3. 
The identity $\varepsilon_W \circ \partial = 0$ then follows 
by considering the boundary of the one-dimensional moduli spaces $\mathcal{M}_W(\gamma)$, 
so that $\varepsilon_W$ is an augmentation.

Given an augmentation $\varepsilon$ of the contact homology DGA $(\mathcal{A}, \partial)$ of a contact manifold, one can define
a linearized complex $(C^\varepsilon, \partial^\varepsilon)$ by taking $C^\varepsilon$ identical to the graded $\mathcal{R}$-module
$C^{\rm cyl}$ from the previous section, and by taking the linearized differential $\partial^\varepsilon : C^\varepsilon \to C^\varepsilon$ 
to be the linear map characterized by
\[
\partial^\varepsilon q_{\gamma^+} = m_{\gamma^+} \sum_{r \ge 0} 
\sum_{\substack{\{\gamma^-_1, \dots, \gamma^-_r\}  \subset \mathcal{P}^g_\alpha \\ |q_{\gamma^-_1}| + \dots + |q_{\gamma^-_r}| 
= |q_{\gamma^+}|-1}} \sum_{k=1}^r
n(\gamma^+;\gamma^-_1, \dots, \gamma^-_r)  
\varepsilon(q_{\gamma^-_1} \dots q_{\gamma^-_{k-1}}) q_{\gamma^-_k} \varepsilon(q_{\gamma^-_{k+1}} \dots q_{\gamma^-_{r}}). 
\]
Since $\varepsilon$ is a map of degree $0$, all terms in which $|q_{\gamma^-_i}| \neq 0$ for some $i \neq k$ vanish in the above expression.
It follows from the properties of $\partial$ and $\varepsilon$ that $\partial^\varepsilon \circ \partial^\varepsilon = 0$. The homology of this
linearized complex is denoted by $CH^\varepsilon(M,\xi)$ and is called the linearized contact homology of $(M, \xi)$ with respect to 
$\varepsilon$. Note that this homology depends on the choice of the augmentation $\varepsilon$ for the contact homology DGA, but that
the collection of all linearized contact homologies for all augmentations of this DGA is an invariant of the contact manifold $(M, \xi)$.

More precisely, following 
\cite[Section~14.5]{arxiv:2307.09068},
given two augmentations $\varepsilon_1$ and $\varepsilon_2$, we define a 
$(\varepsilon_1, \varepsilon_2)$-derivation as a linear map $K : \mathcal{A} \to \mathcal{R}$ such that
\[
K(q_1 \dots q_k) = \sum_{g \in S_k} \sum_{j=1}^k (-1)^{|q_{g(1)} \dots q_{g(j-1)}|} \varepsilon_1(q_{g(1)} \dots q_{g(j-1)})
K(q_{g(j)}) \varepsilon_2(q_{g(j+1)} \dots q_{g(k)}),
\]
for all generators $q_1, \dots, q_k$ of $\mathcal{A}_+$, where $S_k$ denotes the permutation group of $\{ 1, \dots, k \}$.
We then say that $\varepsilon_1$ and $\varepsilon_2$ are DGA-homotopic if there exists a $(\varepsilon_1, \varepsilon_2)$-derivation
$K : \mathcal{A} \to \mathcal{R}$ such that $\varepsilon_1 - \varepsilon_2 = K \circ \partial$. It then follows from homological algebra arguments
that $CH^\varepsilon(M,\xi)$ depends on $\varepsilon$ only through its DGA-homotopy class.

Moreover, if $\Phi : (\mathcal{A}_+, \partial_+) \to (\mathcal{A}_-, \partial_-)$ is a unital DGA map and if $\varepsilon_-$ is an augmentation of 
$(\mathcal{A}_-, \partial_-)$, then its pullback $\varepsilon_+ = \Phi^* \varepsilon_- = \varepsilon_- \circ \Phi$ is an augmentation of 
$(\mathcal{A}_+, \partial_+)$, and its DGA homotopy class depends only on the DGA homotopy class of $\varepsilon_-$ and on the homotopy
class of the map $\Phi$. The invariance properties announced above for the collection of linearized contact homologies then follow from the 
invariance of contact homology discussed in Section~3.1. 
 
If the augmentation is of the form $\varepsilon_{\rm cyl}$ as above, then the linearized differential coincides with the cylindrical differential,
so that linearized contact homology with respect to $\varepsilon_{\mathrm{cyl}}$ is nothing but cylindrical contact homology.  On the other hand,
if the augmentation is of the form $\varepsilon_W$ as above, then arguments formally similar to those described at the end of 
Section~3.1 
show that the DGA-homotopy class of $\varepsilon_W$ depends only on the symplectic filling $(W, \omega)$ and not
on extra choices such as a compatible almost complex structure. In this situation, it is convenient to denote the resulting linearized contact 
homology as $CH(W, \omega)$. A beautiful calculation of this homology, in the case of subcritical Stein manifolds $(W, \omega)$ with its first 
Chern class vanishing on $\partial W$, was obtained by 
Yau
as the main result of
her PhD thesis under Eliashberg's supervision \cite{mr:2700260}.
This case is actually at the intersection of the two special situations
described above, as it turns out that $\varepsilon_W = \varepsilon_{\mathrm{cyl}}$ 
for such manifolds.

\begin{Thm3.1}
[Yau 
\cite{mr:2087083}]
Let $(W, \omega)$ be a subcritical Stein manifold of dimension $2n \ge 4$ such that the restriction to $\partial W$ of its first Chern 
class $c_1(W, \omega)$ vanishes. Then 
\begin{equation} \label{eq:CHStein}
	CH_i(W,\omega) \cong \bigoplus_{m \ge 0} H_{2(n+m-1)-i}(W)
	{\tag{3.3}}
\end{equation}
for all $i \in \mathbb{Z}$.
\end{Thm3.1}

Moreover, He 
\cite{mr:3045340}
showed 
that the contact homology of the boundary of such a subcritical Stein manifold is isomorphic 
to the exterior algebra of the above graded module \eqref{eq:CHStein}. In other words, the terms in the contact homology
differential \eqref{eq:d} corresponding to holomorphic curves with multiple negative punctures have no effect on the resulting 
homology. 

Before this result, contact homology had been computed for a very limited number of examples. 
Fabert
\cite{mr:2563686}
computed the contact homology for mapping tori of Hamiltonian symplectomorphisms of closed, symplectically aspherical manifolds
$(W, \omega)$. Once again, it turns out that only holomorphic cylinders contribute to the contact homology differential, so that 
contact homology is generated as an algebra by countably many copies of the singular homology of $W$.

Yau 
\cite{mr:2597305}
also computed 
the cylindrical contact homology in the case of contact $3$-manifolds supported by an open book decomposition
having as page a punctured $2$-torus and as monodromy a positive Dehn twist along an embedded nonseparating loop.

\subsection*{\textup{3.4.} Relation to symplectic invariants}

Given a symplectically aspherical symplectic manifold $(W, \omega)$ with contact type boundary, one can define its 
linearized contact homology $CH(W,\omega)$. 
It is natural to ask whether this homological invariant is related to Floer type invariants associated to $(W, \omega)$ using a more symplectic 
point of view. Given a time-dependent Hamiltonian function $H : \mathbb{R}/\mathbb{Z} \times W \to \mathbb{R}$, its Hamiltonian vector field is defined by 
$\imath(X_H^\theta) \omega = dH_\theta$ for all $\theta \in \mathbb{R}/\mathbb{Z}$, where $H_\theta  = H |_{\{ \theta \} \times W}$. We then consider the flow 
$(\varphi^\theta_H)_{\theta \in \mathbb{R}/\mathbb{Z}}$ of this 
time-dependent
vector field. The fixed points of the time-one flow $\varphi_H^1$ are in bijective 
correspondence with 
one-periodic orbits of the time-dependent vector field $X_H^\theta$. We require that for any such fixed point $p$ of 
$\varphi^1_H$, the differential of $\varphi^1_H$ at $p$ has no eigenvalue equal to $1$. The Hamiltonian $H$ is then said to be nondegenerate; 
such Hamiltonians form a dense subset of the set of all Hamiltonian functions. We then consider the module $FC(H)$ freely generated by 
all fixed points of $\varphi_H^1$ and graded by the Conley--Zehnder index of the corresponding orbits. 
After choosing a 
one-parameter family of compatible almost complex structures $(J_\theta)_{\theta \in
\mathbb{R}/\mathbb{Z}}$ on $(W, \omega)$,
we define a differential $\partial^{\rm F}$ similarly to the cylindrical
differential in Section~3.2, 
replacing holomorphic cylinders 
with solutions $u : \mathbb{R} \times \mathbb{R}/\mathbb{Z} \to W$ of the Floer equation 
\[
\frac{\partial u}{\partial s} + J_\theta(u) \left( \frac{\partial u}{\partial
\theta} - X_H^\theta(u) \right) = 0,
\]
and asymptotic to $1$-periodic orbits of $X_H^\theta$ as $s \to \pm \infty$. 

The Floer complex is typically defined with coefficients in the Novikov ring $\Lambda_\omega$ of $(W, \omega)$. This is a completion of the group ring 
$\mathbb{Q}[H_2(W, \mathbb{Z})]$, where infinite sums $\sum_A q_A e^A$ are allowed as long as for all $C > 0$ there are finitely many nonvanishing coefficients 
$q_A$ with $\omega(A) \le C$. This is because the Gromov compactness theorem holds below some energy bound, and the energy of a Floer trajectory 
depends on its homology class, not only on its asymptotes as in a symplectization where the symplectic form is exact.

When $W$ is closed, $(FC(H), \partial^{\rm F})$ is called the Floer
complex; its homology is independent of the choices of $H$ and of
$(J_\theta)_{\theta \in \mathbb{R}/\mathbb{Z}}$, and turns out to be isomorphic to 
the singular homology of $W$. When $(W, \omega)$ has contact type 
boundary, we pick a Liouville vector field $V$ near the boundary $M$ so that it is equipped with a contact form $\alpha = \imath(V) \omega$.
We can then complete the symplectic manifold $(W, \omega)$ by attaching to it the positive part of a symplectization 
$(\mathbb{R}^+ \times M, d(e^t \alpha))$ along the boundary of $W$. Then, one has to impose some asymptotic conditions for the geometric objects 
introduced above. First, we require that for $t$ sufficiently large each $J_\theta$ for $\theta \in
\mathbb{R}/\mathbb{Z}$ 
coincide 
with some 
almost complex structure on the symplectization of $(M, \alpha)$ as in
Section~2.2. 
Second, we require that for $t$ sufficiently large
the Hamiltonian $H$ 
be time-independent, that it not depend on the $M$ factor anymore and that it increase quadratically in $t$. In this region, the Hamiltonian 
vector field will be a multiple of the Reeb vector field, with a proportionality constant that increases to infinity with $t$.  Under these conditions,
$(FC(H), \partial^{\rm F})$ is a chain complex. Its homology is called symplectic homology $SH(W,\omega)$ and is independent of the choices of 
$H$ and of $(J_\theta)_{\theta \in \mathbb{R}/\mathbb{Z}}$ with the required properties. As an alternative, for $t$
sufficiently large in $(\mathbb{R}^+ \times M, d(e^t \alpha))$, 
one can also use time-independent Hamiltonians $H$ increasing linearly at infinity, with a slope which is not the period of a closed Reeb orbit. 
Taking a direct limit of the homologies corresponding to an increasing sequence of Hamiltonians with asymptotic slopes going to infinity, one 
recovers the same symplectic homology $SH(W, \omega)$.

There are two essential differences between the constructions of symplectic
homology and of linearized contact homology. First, although 
one-periodic
Hamiltonian orbits coincide with closed Reeb orbits for $t$ sufficiently large in $(\mathbb{R}^+ \times M, d(e^t \alpha))$, there are other types of generators
in the Floer complex. If one chooses a negative $C^2$-small Hamiltonian function in $W$, which depends only on $t$ and grows quadratically in
$(\mathbb{R}^+ \times M, d(e^t \alpha))$, then $FC(H)$ has two types of generators: those corresponding to the Morse complex of $W$ and those 
corresponding to closed Reeb orbits. Moreover, if we assume that $(W, \omega)$
is symplectically atoroidal, or, in other words, that 
$\int_{T^2} f^*\omega = 0$ for all smooth maps $f : T^2 \to W$, then one can
define the action of one-periodic Hamiltonian orbits by 
\[
\mathcal{A}_H(\gamma) = \int_{\mathbb{R} \times S^1} f^*\omega - \int_0^1 H(\theta, \gamma(\theta)) \, d\theta,
\]
where $f$ is a smooth homotopy from
$\gamma$ to a fixed representative of its homotopy class. Then, the action of the former type of orbits is negative, while the action of the latter type of 
orbits is positive. Since the action decreases under the differential $\partial^{\rm F}$, the former type of orbits generates a subcomplex of the
Floer complex, and one can show that its homology is isomorphic to the relative singular homology $H_{*+n}(W, \partial W; \Lambda_\omega)$ with 
Novikov coefficients, with a grading shift by $n = \frac12 \dim W$.
The homology of the corresponding quotient complex is denoted by $SH^+(W, \omega)$ and is called the positive part of symplectic homology.
We then have the following tautological exact sequence:
\begin{equation} 
	\label{eq:SHtauto}
\cdots \to H_{k+n}(W, \Lambda_\omega) \to SH_k(W,\omega) \to SH^+_k(W, \omega)
	\to H_{k-1+n}(W, \Lambda_\omega) \to \cdots.\tag{3.4}
\end{equation}

The second essential difference between symplectic homology and linearized contact homology is that the Hamiltonian vector field is 
time-depen\-dent, and generators of the Floer complex correspond to parametrized orbits, while the Reeb vector field is independent of time
and generators of the contact complex correspond to unparametrized orbits, or,
in other words, to circles of parametrized orbits, as the 
starting point of the parametrization does not matter. 
By comparison with symplectic homology, linearized contact homology can intuitively be thought
of as a quotient
theory, since generators of its chain
complex
are considered modulo the circle action induced by the Reeb flow along periodic
orbits.
Comparing with exact sequences in singular homology in the context
of (semifree) circle actions, it should not be so surprising that we have an exact sequence 
\cite{mr:2471597}
\begin{equation} \label{eq:SHCH}
\cdots \to SH^+_{k-(n-3)}(W, \omega) \to CH_k(W,\omega) \to CH_{k-2}(W, \omega) 
	\to  SH^+_{k-1-(n-3)}(W, \omega) \to \cdots. \tag{3.5}
\end{equation}

In fact, one can formalize this intuition by defining an $S^1$-equivariant version of symplectic homology, denoted by $SH^{S^1}(W, \omega)$, 
where $S^1$ acts by translation on the time variable $\theta \in S^1$. This version of symplectic homology also has a tautological exact sequence
similar to Equation \eqref{eq:SHtauto}, in which the positive part of $S^1$-equivariant symplectic homology is denoted by $SH^{S^1, +}(W, \omega)$. As
in algebraic topology, symplectic homology and its $S^1$-equivariant version satisfy a Gysin-type exact sequence, and in the case of the positive 
parts of these 
homologies, we have an isomorphism of exact sequences, corresponding to the following commutative diagram 
\cite{mr:3671507}:
\[
\xymatrixcolsep{1pc}\xymatrix{
\cdots \ar[r] & SH^+_{k-(n-3)} \ar[r] \ar@{=}[d] & CH_k \ar[r] \ar[d]^{\cong} & CH_{k-2} \ar[r] \ar[d]^{\cong} &  SH^+_{k-1-(n-3)} \ar[r] \ar@{=}[d] & \cdots \\
\cdots \ar[r] & SH^+_{k-(n-3)} \ar[r] & SH^{S^1,+}_{k-(n-3)} \ar[r] & SH^{S^1,+}_{k-2-(n-3)} \ar[r] & SH^+_{k-1-(n-3)} \ar[r] & \cdots 
}
\]
where we removed all decorations $(W, \omega)$ for brevity.

As a consequence of this result, one sees that the positive part of $S^1$-equivariant symplectic homology can be used as a substitute for linearized
contact homology, even when the latter is not well defined 
owing to transversality issues for moduli spaces of holomorphic curves. It is indeed possible 
in a number of situations to find generic almost complex structures and Hamiltonian functions to make the moduli spaces for symplectic homology 
regular, because these geometric structures are time-dependent, unlike in contact homology. In situations where some Reeb vector field for $(M, \xi)$
does not have any contractible orbit, it is possible to define $SH^{S^1, +}(M, \xi)$ using the symplectization instead of a symplectic filling. This invariant 
has mostly been used when $(W, \omega)$ is a Liouville manifold. In this spirit, Gutt and Hutchings 
\cite{mr:3868228}
used this invariant to construct a sequence of
symplectic capacities for star-shaped domains in $\mathbb{R}^{2n}$, leading to
new applications to some symplectic embedding problems. 

Using this relation between linearized contact homology and the positive part of $S^1$-equiva\-riant symplectic homology, it is possible to revisit 
Theorem~3.1 
and to give an alternative proof via symplectic geometry. The key point is that symplectic homology, 
as well as its 
$S^1$-equivariant counterpart, both vanish in the case of a subcritical Stein
manifold. This is 
due to the fact that these manifolds factor as the
product of a Stein manifold with $\mathbb{C}$, as shown by 
Cieliebak
\cite{arxiv:math/0204351};
see also 
\cite[Section~14.4]{mr:3012475}.
The symplectic topology of
Stein manifolds is also the 
subject of another chapter of this volume
\cite{cm:cieliebak.2024}.
The conclusion for symplectic homology then follows from the
K\"unneth formula proved by 
Oancea
\cite{mr:2208949}
and the fact that symplectic homology of $\mathbb{C}$ vanishes. The Gysin exact sequence between
symplectic homology and its $S^1$-equivariant version then shows that the latter has the same rank in any degree of a given parity. But an easy
inspection shows that the generators of its chain complex have nonnegative grading, so that $S^1$-equivariant symplectic homology vanishes as well.
Once this is established, the tautological exact sequence \eqref{eq:SHtauto} in the $S^1$-equivariant case shows that the positive part of 
$S^1$-equivariant symplectic homology (or in other words linearized contact homology) is isomorphic to $H^{S^1}_*(W)$. Since the
circle acts on the $\theta \in S^1$ variable of the Hamiltonian and not on $W$, the latter homology is isomorphic to 
$H_*(W) \otimes H_*(BS^1)$, and this coincides with the result given in
Theorem~3.1. 

\subsection*{\textup{3.5.} Variants of contact homology}
In order to complete this brief tour of the foundations of contact homology, let us mention some variants of the above constructions that led to other 
interesting homological invariants in contact and symplectic topology.

Embedded Contact Homology (ECH) was introduced by 
Hutchings;
for an introduction
see, for example, 
\cite{mr:3220947}.
This theory also consists in a chain complex whose 
differential counts holomorphic curves in a symplectization, but as suggested by the name ECH, only embedded holomorphic curves are considered.
It is also a four-dimensional theory, because embedding control is made possible
by positivity of intersection, which is a four-dimensional phenomenon.
Another major difference between ECH and contact homology is that the former turns out to be a topological invariant of the $3$-manifold under 
consideration, and does not depend on the contact structure which is used to
define it, because it is isomorphic to the Seiberg--Witten Floer cohomology
defined by 
Kronheimer
and 
Mrowka
\cite{mr:2388043}.
This isomorphism was established by 
Taubes,
who used it to prove 
\cite{mr:2350473}
the Weinstein conjecture in 
dimension three. This very important conjecture will be discussed in more details
in Section~4.4. 
ECH is also a formidable tool to study 
symplectic embedding problems and to obtain sharp restrictions on them via the so-called ECH capacities. All this is well outside our present scope and 
certainly deserves its own volume.

Sutured contact homology was introduced by 
Colin,
Ghiggini,
Honda and Hutchings 
\cite{mr:2851076}.
It is a generalization of contact homology 
which is defined for a certain class of contact manifolds with boundary, called sutured contact manifolds. The latter are the contact analogue of 
symplectic manifolds with contact type boundary. Given a contact manifold $(M, \xi)$ with boundary, one requires that the $\partial M$ 
be 
transverse to a vector field preserving $\xi$. Then, one restricts to contact forms $\alpha$ for $\xi$ so that the tangency locus of the Reeb 
field $R_\alpha$ to $\partial M$ is a smooth submanifold $\Gamma$ on which $\alpha$ restricts to a contact form: it is the suture of $(M, \xi)$.
Its complement $\partial M \setminus \Gamma$ consists of regions $R_\pm$ where $\pm R_\alpha$ points outside $M$, and $(R_\pm, \alpha)$ 
are Liouville manifolds with convex boundary $\Gamma$. Finally, the Reeb field $R_\alpha$ is required, along the suture $\Gamma$, to point 
from $R_-$ to $R_+$ ; this condition corresponds to a convexity condition for the contact form $\alpha$ along $\partial M$. One can then complete
such a contact manifold $(M, \alpha)$ so that all closed Reeb orbits lie in the original manifold and all holomorphic curves in the corresponding 
symplectization that are asymptotic to these Reeb orbits do not leave a compact subset of the completed manifold. Sutured contact homology can
then be defined in a similar way to contact homology. 
Golovko
\cite{mr:3320870}
computed the sutured contact homology for universally tight contact structures 
on solid tori. This variant of contact homology can also lead to Legendrian invariants by considering the sutured contact homology of the complement
of Legendrian submanifolds.

Local contact homology was introduced by 
Hryniewicz
and 
Macarini
\cite{mr:3326300}.
This is a homological invariant associated to a single closed Reeb orbit 
$\gamma$, which is possibly degenerate. Intuitively, it measures the local contribution of this closed Reeb orbit to cylindrical or linearized contact homology 
of the whole contact manifold, after a perturbation of the contact form in order to replace $\gamma$ with a collection of nearby nondegenerate closed 
Reeb orbits. It turns out that the rank of this invariant is bounded in any given grading and for any iteration of the orbit $\gamma$. Applications of this
invariant include a generalization of Gromoll--Meyer’s theorem on the existence of infinitely many simple periodic orbits, resonance relations and 
conditions for the existence of nonhyperbolic periodic orbits.

Momin
\cite{zbl:1248.37057}
defined a version of cylindrical contact homology for orbit complements in
dimension three. More precisely, given a link $L$ in a contact 
$3$-manifold $(M, \xi)$, we 
restrict 
ourselves to contact forms $\alpha$ such that $L$ is realized by a collection of closed Reeb orbits. We define a graded 
module that is linearly generated by closed Reeb orbits $\gamma$ of such contact forms in
$M \setminus L$, such 
that $\gamma$ realizes 
a fixed free homotopy class $[a]$ in this manifold.
It is equipped with differential counting rigid holomorphic cylinders that are contained in the symplectization of $M \setminus L$. Under suitable assumptions 
on $L$, $\alpha$ and $[a]$, this is a chain complex whose homology has some invariance properties. This invariant can be used to establish some implied
existence results: assuming that the Reeb flow admits some closed orbits with certain properties, one deduces that it has a certain number of other 
closed orbits in some free homotopy class.

The Euler characteristic is a convenient numerical invariant for chain complexes of finite rank. When the chain complex for linearized or cylindrical 
contact homology has finite rank in any given degree, but has a total rank which is infinite, it is sometimes possible to define its mean Euler characteristic.
In the case of cylindrical contact homology, this was defined by van~Koert 
\cite{cm:vKthesis.2005} 
as
\[
\chi(M,\xi) = \lim_{N \to \infty} \frac1N \sum_{m = -N}^N (-1)^m \textrm{rank } CH^{\rm cyl}_m(M,\xi).
\]
A huge advantage of this invariant is that it can be computed in terms of closed
Reeb orbits and their Conley--Zehnder indices, without having to compute
the contact homology differential, or, in other words, 
having to find all rigid holomorphic curves. Going beyond distinguishing contact
manifolds,
in 
\cite{mr:2998911}
Frauenfelder,
Schlenk
and 
van~Koert
used
this invariant 
to study some contact embedding problems.
This work was then used by 
Kang
\cite{mr:3158574}
to prove the existence of at least two closed Reeb orbits on some classes of
contact manifolds (see Section~4.4 
below for more advanced 
results in this direction). 
The mean Euler characteristic was also used 
by 
Chiang,
Ding
and 
van~Koert 
\cite{mr:3210581}
to prove 
that some fibered Dehn twists are not symplectically isotopic to the identity
relative to the boundary.
Espina
\cite{mr:3215222}
extended this notion to a broader class of contact manifolds; she also studied the effect of subcritical contact surgery on this invariant.

\section*{4. Applications of contact homology}

\subsection*{\textup{4.1.} Distinguishing contact structures}

Since contact homology and its variants are invariants of contact manifolds, their first applications consist in distinguishing contact structures on a 
given closed manifold $M$. The first application of this nature was obtained by 
Ustilovsky
as the main result of his PhD thesis under 
Eliashberg's supervision 
\cite{mr:2700259},
with $M = S^{4k+1}$ for any integer $k \ge 1$. The  Brieskorn manifolds, defined as the transverse intersection of
the singular complex hypersurface $z_0^{a_0} + z_1^{a_1} + \dots + z_n^{a_n} = 0$ with a small unit sphere around the origin, are diffeomorphic
to $S^{4k+1}$ when $n = 2k+1$, $a_1 = \dots = a_{2k+1} = 2$ and $2 < a_0 = p \equiv \pm 1 \mod 8$. The $1$-forms 
\begin{equation} \label{eq:brieskorn}
\alpha_{(a_0, \dots, a_n)} = \frac{i}8 \sum_{j=0}^n a_j (z_j d{\bar z}_j - {\bar
	z}_j dz_j)\tag{4.1}
\end{equation}
are contact forms and we denote by $\xi_p$ the corresponding contact structures on $S^{4k+1}$ in the above particular case.

Note that contact structures can sometimes be distinguished by their classical invariants, defined using tools from algebraic topology. 
These are all determined by the formal homotopy class of the contact structure, which is defined as its homotopy class as a complex vector
subbundle of the tangent bundle of the 
ambient 
manifold. For spheres of dimension $4k+3$ with $k \ge 1$, there are infinitely many formal 
homotopy classes, and one can show that the contact structures obtained as above
fall into distinct classes, so 
contact homology is not 
needed to distinguish those contact structures. On the other hand, for spheres of dimension $4k+1$ with $k \ge 1$, there are only finitely many 
such formal homotopy classes, so the above construction leads to infinitely many contact structures having the same formal homotopy classes. 
In this case, it is very useful to compute their cylindrical contact homology.

\begin{Thm4.1}
[Ustilovsky %
\cite{mr:1704176}]
For all $k \ge 1$ and $p \ge 2$ such that $p \equiv \pm 1 \mod 8$,
there exists a contact form $\alpha_p$ for $(S^{4k+1}, \xi_p)$ such that cylindrical contact homology over $\mathcal{R} = \mathbb{Q}$ 
is well defined and the dimension $d_m$ of $CH_m^{\rm cyl}(S^{4k+1}, \alpha_p)$ is given by
\[
d_m = \left\{
\begin{array}{ll}
0 & \textrm{if } m \textrm{ is odd or } m < 2n-4, \\
2 & \textrm{if } m = 2 \lfloor 2N/p \rfloor + 2(N+1)(n-2) \\
& \quad \textrm{for some integer } N \ge 1 \textrm{ such that } 2N+1 \notin p\mathbb{Z}, \\
1 & \textrm{otherwise.}
\end{array} \right.
\]
\end{Thm4.1}

This explicit calculation was possible thanks to the fact that the Conley--Zehnder indices of all closed Reeb orbits of $\alpha_p$, which is a
suitable perturbation of the $1$-form $\alpha_{(a_0, \dots, a_n)}$ in 
Equation \eqref{eq:brieskorn}, have the same parity, 
so that $\partial^{\rm cyl} = 0$. 
The result of this calculation shows in particular that cylindrical contact homology remembers the value of $p$, and the invariance of contact homology
then implies that $S^{4k+1}$ admits infinitely many pairwise nondiffeomorphic contact structures.

At the time Ustilovsky's paper 
\cite{mr:1704176}
appeared, 
invariance of cylindrical contact homology in the presence of contractible closed Reeb orbits was not fully established, due to
the possible presence of multiply covered holomorphic planes with a positive asymptote at such orbits, for which transversality cannot be achieved
with a generic almost complex structure. Using the considerations explained in
Section~3.4, 
Gutt
\cite{mr:3734608}
gave a complete proof of 
the existence of infinitely many pairwise nondiffeomorphic contact structures on $S^{4k+1}$ using the positive part of $S^1$-equivariant symplectic 
homology as a well-established substitute for cylindrical contact homology in
Theorem~4.1. 

Using the same techniques, 
Fauck
\cite{mr:4192450}
extended this result to all contact manifolds $(M, \xi)$ of dimension $4k+1$ admitting an exact 
symplectic filling such that noncontractible loops in $M$ are not contractible in the filling either, and satisfying an additional assumption 
(asymptotically finitely generated) on the quantity of closed Reeb orbits with a
given Conley--Zehnder index. 

In dimension three, one can construct infinitely many nondiffeomorphic contact structures on a toroidal manifold $M$ by cutting this manifold
along an incompressible torus and inserting a domain of the form $T^2 \times [0, 2k\pi]$ for some integer $k \ge 1$ equipped with the contact structure 
$\xi_k = \ker (\cos(\theta) \,dx + \sin(\theta) \,dy)$, where $x$ and $y$ are global angular coordinates on $T^2$ and $\theta \in [0, 2k\pi]$. These contact
structures can then be distinguished by their Giroux torsion, which is defined as the largest integer $k$ such that the above contact domain
can be embedded in a given contact manifold. As this invariant is difficult to compute, it is simpler to use cylindrical contact homology, and this 
was done 
in 
\cite{mr:2116317}
by Colin and the present author. 
In dimensions greater than three, 
Massot,
Niederkr\"uger
and 
Wendl 
\cite{mr:3044125}
generalized this construction
by producing contact structures $\xi_k$ on $T^2 \times M$ for any integer $k \ge 1$, where $M$ are suitable contact manifolds of arbitrary dimensions, 
equipped with 
two contact forms $\alpha_+$ and $\alpha_-$, in such a way that 
\[
	\xi_k = \ker \Bigl(\frac{1+ \cos kx}2 \alpha_+ + \frac{1-\cos kx}2 \alpha_- 
+ \sin kx \, dy\Bigr). 
\]
As in dimension three, these contact structures are distinguished using cylindrical contact homology. 

In order to compute 
the contact homology or one of its variants, it is useful to
take advantage of the symmetries of natural contact forms and their associated Reeb 
vector fields. For instance, the Reeb flow associated to the contact form
$\alpha_{(a_0, \dots, a_n)}$ in Equation \eqref{eq:brieskorn} is periodic, though all 
simple closed Reeb orbits do not have the same period. In other words, the quotient of $S^{4k+1}$ by the Reeb flow is a symplectic orbifold. In more
general situations, one could have submanifolds of $(M,\xi)$ that are foliated by closed Reeb orbits of different periods. Such contact forms are 
degenerate, but one can impose 
the condition that the eigenspace of the eigenvalue $1$ coincide with the tangent space of these foliated submanifolds, in order to
retain a nondegeneracy property in the transverse directions. This is analogous
to the notion of a Morse--Bott function. It is possible to compute 
contact homology using such a contact form 
\cite{mr:1969267}.
To summarize, the generators of the Morse--Bott complex for contact homology are critical points of 
auxiliary Morse functions chosen on the foliated submanifolds.  The differential 
counts cascades, which consist of an 
alternance between fragments of 
gradient trajectories for the auxiliary Morse functions and $J$-holomorphic curves for an almost complex structure $J$ which is compatible with the 
Morse--Bott geometric data. The fragments joining two holomorphic curves can have an arbitrary finite length, while those reaching the generators 
of the complex are parametrized by a half-line. For more details, see, for
example,
\cite{mr:2555933}.

In particular, the calculation of cylindrical contact homology in Theorem~4.1 
can be made using the contact form $\alpha_{(a_0, \dots, a_n)}$ 
in Equation \eqref{eq:brieskorn} without any perturbation. The corresponding Reeb flow is given by $\varphi^{R_\alpha}_t(z_0, \dots, z_n) = (e^{4it/p} z_0,
e^{2it} z_1, \dots, e^{2it} z_n)$, so that orbits with $z_0 =0$ have minimal period $\pi$ while orbits with $z_0\neq 0$ have minimal period $p\pi$. 
It turns out that the space of closed Reeb orbits with period $kp\pi$ is homeomorphic to $\mathbb{C}\mathrm{P}^{n-1}$ for all $k \ge 1$, while the space of closed
Reeb orbits with period $k \pi$ has the same singular homology as $\mathbb{C}\mathrm{P}^{n-2}$ for all $k \ge 1$ not multiple of $p$. Due to the $S^1$-symmetry 
induced by the Reeb flow, there are no rigid holomorphic curves, so that cylindrical contact homology will be a repetition of the singular homology 
of the orbit spaces for periods $(kp+1)\pi, \dots (kp+p-1)\pi, (k+1)p\pi$ for $k
\ge 0$, with suitable grading shifts due to the Conley--Zehnder index.
It follows easily that there is a repetitive pattern with respect to the grading for the dimensions of cylindrical contact homology, a fact which is not
completely obvious from the statement of Theorem~4.1 

More generally, van~Koert 
\cite{mr:2395966}
computed the cylindrical contact homology for Brieskorn manifolds with many different values 
of the exponents $a_0, \dots, a_n$. As a consequence of this, he exhibited large classes of contact manifolds admitting infinitely many 
nondiffeomorphic contact structures. More results on this class of contact manifolds were surveyed by 
Kwon
and van~Koert 
\cite{mr:3483060},
using $S^1$-equivariant symplectic homology. Similar computations were used by 
Boyer,
Macarini and van~Koert 
\cite{mr:3543703}
to study the
space of positive Sasakian structures on these manifolds.

Other applications of the above Morse--Bott techniques were obtained by Abreu and Macarini 
\cite{mr:2881318}
who computed 
the
cylindrical contact homology combinatorially for large classes of toric contact manifolds. In particular they exhibited infinitely many 
nondiffeomorphic contact structures on $S^2 \times S^3$ with vanishing first Chern class. In a similar spirit, Boyer and 
Pati
\cite{mr:3207586}
established criteria for contact structures on $S^3$-bundles over $S^2$ to be nondiffeomorphic, even when they are not distinguished 
by their first Chern class. Returning to toric contact manifolds, 
Marinkovi\'c
\cite{mr:3729335}
related the cylindrical contact homologies of two such
manifolds, when one is obtained from the other by a contact 
blow-up.

\subsection*{\textup{4.2.} Topology of the space of contact structures}

In view of  
Gray's stability theorem, distinguishing contact structures amounts to the study of the set of connected components of the space of contact 
structures, modulo the contact mapping class group. It turns out that contact homology can also be used to study higher homotopy groups of 
the space of contact structures. Given an odd dimensional manifold $M$, Let us denote by $\Xi(M)$ the set of all contact structures on $M$,
equipped with the topology induced by the Grassmannian manifold of all codimension one subbundles of the tangent bundle $TM$. Given a
contact structure $\xi_0$ on $M$, one can think of elements in $\Pi_k(\Xi(M), \xi_0)$ as homotopy classes of families $\xi_{(t_1, \dots, t_k)}$
of contact structures on $M$ parametrized by $[0,1]^k$, 
and which coincide with $\xi_0$ if $t_j = 0$ or $1$ for some $j \in \{ 1, \dots, k \}$. 
Given such a $k$-parameter family, one can construct a $k-1$-parameter family of symplectic cobordisms
$(\mathbb{R} \times M, \omega_{(t_2, \dots, t_k)})$, 
all interpolating from $(M, \xi_0)$ to itself, and which coincide with the symplectization of $(M, \xi_0)$ if $t_j = 0$ or $1$ for some $j \in \{ 2, \dots, k \}$. 
To carry out this construction, 
pick a family of contact forms $\alpha_{(t_1, \dots, t_k)}$ for $\xi_{(t_1, \dots, t_k)}$ which coincides with a fixed contact form $\alpha_0$
for $\xi_0$ when $t_j = 0$ or $1$ for some $j \in \{ 1, \dots, k \}$. Choose a smooth decreasing function $f :
\mathbb{R} \to [0,1]$ such that $f(t) = 0$ 
if $t$ is sufficiently large, $f(t) = 1$ if $t$ is sufficiently small and $f$ varies very slowly. Then 
$\omega_{(t_2, \dots, t_k)} = d(e^t \alpha_{(f(t), t_2, \dots, t_k)})$ is the desired family of symplectic forms.

Counting rigid holomorphic curves in this family of symplectic cobordisms leads to maps descending to homotopy groups of $\Xi(M)$. 
More precisely, one obtains group morphisms 
\cite{mr:2200047}
\[
\eta_1 : \pi_1(\Xi(M), \xi_0) \to \textrm{Aut}(CH^{\rm cyl}(M, \xi_0))
\]
and
\[
\eta_k : \pi_k(\Xi(M), \xi_0) \to \textrm{End}_{1-k}(CH^{\rm cyl}(M, \xi_0))
\]
for $k \ge 2$, where for a graded module $V$, $\textrm{Aut}(V)$ denotes the nonabelian group of degree $0$ automorphisms of $V$ and
$\textrm{End}_{1-k}(V)$ denotes the additive group of degree $1-k$ endomorphisms of $V$. These morphisms can be used to detect nontrivial
elements in the homotopy groups of $\Xi(M)$. For example, $\pi_1(\Xi(T^3), \xi_0)$ contains an infinite cyclic subgroup for any tight contact 
structure $\xi_0$. This subgroup is generated by loops $t \in \mathbb{R}/2\pi\mathbb{Z} \mapsto \ker(\cos(nz-t) \,dx + \sin(nz-t) \,dy) \in \Xi(T^3)$, where $(x,y,z)$
are angular coordinates on $T^3$ and $n$ is a positive integer. A similar statement holds for $T^2$-bundles of $S^1$ and for $T^5$, using
appropriate base points $\xi_0$. For a compact, orientable $2n$-manifold $Q$,
$\,\pi_{2n-1}(\Xi(ST^*Q), \xi_0)$ also contains an infinite cyclic 
subgroup, when $\xi_0$ is the canonical contact structure on the unit cotangent bundle $ST^*Q$.

\subsection*{\textup{4.3.} Tight vs. overtwisted}
\label{sec:OT}

In the case of 
three-dimensional contact manifolds $(M^3, \xi)$, Eliashberg introduced the notion of overtwisted contact structures:
a contact structure $\xi$ is said to be overtwisted if there exists a smooth embedding of the $2$-disk in $M^3$ such that, along the
boundary of this disk, its tangent space coincides with $\xi$; in that case, this embedded $2$-disk is called an overtwisted disk. 
If no such disk exists in $(M, \xi)$, then $\xi$ is said to be tight. This dichotomy between tight and overtwisted contact structures is 
of special importance because Eliashberg proved 
\cite{mr:1022310}
that overtwisted contact structures satisfy an h-principle 
and are therefore classified in purely topological terms. In contrast, tight contact structures tend to have a geometric interest: 
for example, 
they include, but are not limited to, those contact structures that are induced by a symplectic filling. Since the existence of overtwisted 
contact structures is a purely topological phenomenon, it makes sense that contact homology has nothing to say about them. The following 
result was proved in 
\cite{mr:2230587}
by Yau, while a different proof of the same 
was given in the appendix by Eliashberg.

\begin{Thm4.2}
[Yau, Eliashberg %
\cite{mr:2230587}]
If $(M, \xi)$ is a closed overtwisted contact $3$-manifold, then its contact homology $CH(M, \xi)$ vanishes.
\end{Thm4.2}

In dimension greater than three, there was no notion of overtwisted contact structures at the time of the above result. However, some classes 
of contact manifolds in arbitrary dimensions were singled out as sharing some geometric properties with overtwisted contact manifolds, and
contact homology was proved to vanish for some of these classes. This includes the class of contact structures obtained by a negative 
stabilization of an open book decomposition supporting another contact structure 
\cite{mr:2646902}.
Another significant class consisted of the 
so-called PS-overtwisted contact structures 
\cite{mr:2286033},
defined using a higher dimensional generalization of the overtwisted disk. 
An unpublished result (see 
\cite[Theorem~4.10]{mr:2555933})
showed that contact homology vanishes for such contact manifolds as well.
Moreover, the vanishing of contact homology was shown to be equivalent to the vanishing of the rational SFT or of the full SFT 
\cite{mr:2606237}.

This situation changed drastically when 
Borman,
Eliashberg and 
Murphy
\cite{mr:3455235}
defined the notion of overtwisted contact structure 
in arbitrary dimension and established an h-principle for this class of contact
structures. This celebrated result is the 
subject of another chapter in this volume
\cite{cm:etnyre.2024}. 
Let us simply mention here that subsequent results by 
Casals,
Murphy and 
Presas
\cite{mr:3904160}
implied that contact homology vanishes 
for overtwisted contact manifolds of arbitrary dimensions, and this superseded all of the above-mentioned vanishing results.

In view of these vanishing results, it is natural to ask whether contact homology detects overtwisted contact structures; in other words, does
$CH(M, \xi) = 0$ imply that $\xi$ is overtwisted? This important question was answered negatively by 
Avdek
\cite{mr:4599310}.
More precisely, he showed 
that contact $+1$-surgery on the standard contact $3$-sphere along a Legendrian link having a component which is a right-handed trefoil produces
contact manifolds with vanishing contact homology. One of these contact
manifolds was already known to be tight, 
based on arguments from 
Heegaard Floer homology.

As a consequence of the above groundbreaking result by Borman, Eliashberg and Murphy, it became interesting to construct and study tight contact 
manifolds in dimension greater than three. In particular, a former construction 
\cite{mr:1912277}
of a contact structure on $T^2 \times M$ based on an 
open book decomposition supporting a contact structure on $M$ was studied in this new light by several groups of researchers. 
Bowden, 
Gironella,
Moreno
and 
Zhou
\cite{arxiv:2211.03680}
used contact homology to show that some of these contact manifolds are tight.
This, in turn,
allowed them to construct many examples of tight contact manifolds without
symplectic fillings. More recently, Avdek and Zhou 
\cite{arxiv:2404.16311}
showed
that all contact structures obtained from the above construction are tight, by computing the contact homology of a suitable covering of $T^2 \times M$.

\subsection*{\textup{4.4.} Dynamical complexity}
\label{sec:dynamics}

Reeb vector fields constitute a remarkable class of nonvanishing vector fields, and their study from a dynamical point of view is of great interest. 
In particular, they generalize the geodesic flow in the case of unit cotangent bundles equipped with their canonical contact structures. A central
question about the dynamical properties of Reeb vector fields is the celebrated Weinstein conjecture, which states in its modern form that any
Reeb vector field on a closed contact manifold admits a periodic orbit. This conjecture motivated a lot of intensive work in contact geometry,
leading to a plethora of interesting results. Indeed, contact homology can be
thought of as an especially appropriate tool to prove this conjecture
for a given closed contact manifold $(M, \xi)$. It suffices 
to show that $CH(M,\xi) \neq \mathcal{R}$, or that $CH^{\rm cyl}(M, \xi) \neq 0$,
or that $CH^\varepsilon(M,\xi) \neq 0$ for some augmentation $\varepsilon$, to deduce that the Weinstein conjecture holds for $(M,\xi)$. 
Hofer's proof 
\cite{mr:1244912}
of the Weinstein conjecture for overtwisted contact $3$-manifolds can be reinterpreted as the main step of the proof outlined
in Section~4.3 
that contact homology vanishes for such manifolds. All computations of contact homology and its variants made so far
establish the Weinstein conjecture for the corresponding contact manifolds. However, the case of general closed contact manifolds still eludes us, 
and in comparison with ECH it seems that we are still lacking a general structural result for contact homology in order to prove the Weinstein
conjecture in dimensions greater than 
three.

Instead of looking for a single closed Reeb orbit in very general contact manifolds, one can turn to a more quantitative version of the Weinstein 
conjecture for some specific classes of contact manifolds. The goal here is to establish some (possibly sharp) lower bound for the number of
geometrically distinct closed Reeb orbits on a given contact manifold. Contact homology is a useful tool in that respect, but its rank does not 
directly provide the desired lower bound, as multiple covers of a given closed orbit are considered as separate generators of the contact chain 
complex, but are not geometrically distinct from the underlying simple closed Reeb orbit. In order to extract information about those simple orbits, 
it is necessary to understand 
thoroughly how the Conley--Zehnder index behaves for iterated orbits. The necessary information is contained in the
so-called common index jump theorem, due to Long and Zhu 
\cite[Theorem~4.3]{mr:1906590}.
Combining cylindrical contact homology with this result,
 Abreu and Macarini 
 \cite{mr:3625067},
 Gutt and Kang 
 \cite{mr:3580178},
and 
Ginzburg
and 
G\"urel
 \cite{mr:4100023}
 (using more symplectic techniques),
 proved independently that a 
 star-shaped hypersurface in $\mathbb{R}^{2n}$, equipped with the restriction of the standard contact form, has at least $n$
 simple closed Reeb orbits, under the assumption that the contact form is nondegenerate and dynamically convex. This last assumption means
 that all closed Reeb orbits have Conley--Zehnder index at least $n+1$. Similar
 results hold for prequantization bundles 
\cite{mr:3625067},
 with as lower 
 bound the sum of the Betti numbers of the base, under some assumptions of a similar nature.
 
Another approach to the abundance of closed Reeb orbits is to consider the number $N_T(\alpha)$ of closed Reeb orbits  with period at most $T$
corresponding to the contact form $\alpha$, and to study the growth rate of $N_T(\alpha)$ as $T$ tends to infinity. Since the period of a closed Reeb
orbit coincides with its action, which decreases under the differential, the module generated by $\mathcal{P}^{\le T}_\alpha$ forms a 
subcomplex for contact homology, and the rank of its homology $CH^{\mathrm{cyl}}_{\le T}(M, \alpha)$ provides a lower bound for $N_T(\alpha)$. 
Vaugon
\cite{mr:3342671}
showed that this growth rate is independent of the choice of a contact form for a given contact structure. She also showed 
that this growth rate is exponential in the case of a closed $3$-manifold such that its 
JSJ decomposition contains a hyperbolic component 
that fibers over the circle. This was motivated by work of Colin and Honda 
\cite{mr:3017043},
in which they computed the cylindrical or linearized contact homology 
for contact $3$-manifolds supported by an open book decomposition with special types of monodromies, including periodic or pseudo-Anosov maps. 
They also showed that linearized contact homology has exponential growth for some pseudo-Anosov monodromies.
 
The complexity of the flow $\varphi$ of a vector field, seen as a dynamical
system, can be measured by its topological entropy $h_{\mathrm{top}}(\varphi)$.
It is especially important to be able to distinguish between dynamical systems
having vanishing 
versus positive topological entropy. Using cylindrical 
contact homology, 
Alves
\cite{mr:3590356}
established a positive lower bound for the topological entropy of the Reeb flow in the case of the same class of 
closed $3$-manifolds as above, as well as for contact $3$-manifolds obtained by
a 
Foulon--Hasselblatt surgery on the unit cotangent bundle of a closed 
surface along the Legendrian lift of a separating geodesic. He also proved 
a similar result for closed contact $3$-manifolds admitting a 
Reeb flow which is Anosov 
\cite{mr:3570995}.
Note that other results of the same nature but for other classes of contact manifolds can be obtained using a variant of 
Legendrian contact homology, but the latter tool is outside our present scope. In the same period, 
Foulon,
Hasselblatt
and 
Vaugon
\cite{mr:4353959}
improved the above results to the case of contact $3$-manifolds obtained by a
Foulon--Hasselblatt surgery on the unit cotangent bundle of a closed 
surface along the Legendrian lift of any closed geodesic.

\subsection*{\textup{4.5.} Other applications}

Let us finally mention some other applications of contact homology that we will not develop in detail in this text.

The celebrated Gromov nonsqueezing theorem is one of the cornerstones of symplectic topology. Eliashberg, 
Kim
and 
Polterovich
\cite{mr:2284048}
used contact homology to establish a contact geometric version of this theorem, with a
crucial distinction between 
large-scale and small-scale phenomena. They also discover a relation to the existence 
of a partial order on the universal covering of the contactomorphism group of the contact manifold under consideration.
This initiated a whole new avenue of research, which we do not develop 
further, as it is the subject of 
another chapter of this volume
\cite{cm:uljarevic.2024}.

As we saw in Section~3.4, 
symplectic homology and its $S^1$-equivariant counterpart both vanish in the case of subcritical Stein manifolds.
This also means that the attachment of a subcritical handle to a symplectic manifold with convex boundary does not change these invariants,
and that its effect on linearized contact homology is purely topological. However, the situation drastically changes when attaching a critical handle
along a Legendrian sphere $\Lambda$ in the contact boundary. 
Together with the present author, Ekholm and Eliashberg 
described
the effect of such an operation on the above invariants 
using long exact sequences involving suitable variants of the Legendrian contact homology of $\Lambda$ 
\cite{mr:2916289}.
This first result has many ramifications and applications, which we do not detail here, as it is the 
subject of another chapter of this volume \cite{cm:ekholm.2024}.

\section*{Acknowledgements} 
The author is first and foremost deeply indebted to Yasha Eliashberg for patiently teaching him so many facets of
the wonderful field of contact and symplectic topology. The author is also grateful to all his current and former
doctoral students and postdocs, for their numerous questions and discussions,
and for particular suggestions regarding
explanations to use in this text. Many thanks for the very warm hospitality of the
Laboratoire de Math\'ematiques Jean Leray at Nantes Universit\'e, where most of the writing of this text took
place.  This text could not have been completed without the unending patience of the editors and the loving
support of the author's wife.  

\printbibliography

 \end{document}